\documentclass[11pt]{amsart}

\usepackage{a4wide, amsmath, amsfonts, amssymb, mathrsfs, amsthm, bbm, stmaryrd}
\usepackage[hyperfootnotes=false]{hyperref}
\usepackage{wasysym}

\allowdisplaybreaks[2]

\numberwithin{equation}{section}

\newtheorem{theorem}{Theorem}[section]
\newtheorem{lemma}[theorem]{Lemma}
\newtheorem{corollary}[theorem]{Corollary}
\newtheorem{proposition}[theorem]{Proposition}

\theoremstyle{remark}
\newtheorem{remark}[theorem]{Remark}

 \global\long\def\sbr#1{\left[ #1\right] }
 \global\long\def\cbr#1{\left\{  #1\right\}  }
 \global\long\def\rbr#1{\left(#1\right)}

 \global\long\def\E{\mathbb{E}}
 \global\long\def\P{\mathbb{P}}
 \global\long\def\R{\mathbb{R}}

 \global\long\def\dd#1{\textnormal{d}#1}

 \global\long\def\ra{\rightarrow}
 
 \global\long\def\ns{\infty}

\title[Estimates of moments and tails
of hitting times of Bessel processes]{Moments and tails
of hitting times of Bessel processes and convolutions of elementary mixtures of exponential distributions}

\author[Bednorz]{Witold M. Bednorz}
\address{Witold M. Bednorz, University of Warsaw, Poland}
\email{wbednorz@mimuw.edu.pl}

\author[\L ochowski]{Rafa\l{} M. {}\L ochowski}
\address{Rafa\l\ M. {\L}ochowski, Warsaw School of Economics, Poland}
\email{rlocho@sgh.waw.pl}

\date{\today.}

\begin{document}

\begin{abstract}
We present explicit estimates of right and left tails and exact (up to
universal, multiplicative constants) estimates of tails and moments of hitting
times of Bessel processes. The latter estimates are obtained from more general estimates of moments and tails established for convolutions of elementary mixtures of exponential distributions.
\end{abstract}

\maketitle

\noindent\textbf{Keywords:} Bessel process, hitting times, moments, tails. \\
\textbf{MSC 2020 Classification:} 60G40, 60J60.


\section{Introduction} \label{sec:Introduction}

Let $X_{t}$, $t\ge0$ be a $\delta$-dimensional Bessel process or
a Bessel process with index $\nu=\delta/2-1>-1$, starting from $x_{0}\ge0$,
that is nonnegative process satisfying $X_{t}^{2}=Z_{t}$ where $Z_{t}$,
$t\ge0$ is a process satisfying the following stochastic differential
equation driven by a standard, one dimensional Wiener process $\beta_{s}$,
$s\ge0$:
\begin{equation}
Z_{t}=x_{0}^{2}+2\int_{0}^{t}\sqrt{\left|Z_{s}\right|}\dd{\beta_{s}}+\delta\cdot t.\label{eq:Bessel}
\end{equation}
For any $z_{0}=x_{0}^{2}\ge0$ and $\delta\ge0$ equation (\ref{eq:Bessel})
has unique strong solution, which is non-negative, cf. \cite[Sect. 1, Chapt. XI]{RevuzYor:2005}.

Let $c$, $\delta$ be fixed positive reals and $z_{0}\in[0,c^{2})$.
In this paper we investigate the laws of hitting times 
\begin{equation}
\tau^{c}:=\inf\cbr{ t\ge0:Z_{t}=c^{2} } =\inf \cbr{ t\ge0:X_{t}=c } .\label{eq:tau_c_def}
\end{equation}
More precisely, we investigate concentration of measure phenomenon
of these laws around the mean of $\tau^{c}$ which reads as 
$
\E\tau^{c}=\rbr{c^{2}-z_{0}}/{\delta}.
$

Investigation of tails of hitting time $\tau^{c}$ has a long history.
Ciesielski and Taylor \cite{CiesTaylor:1962} found an exact formula
for the tails of $\tau^{c}$ in the case when $\delta$ is a positive
integer and $z_{0}=0$ but the problem was already investigated by
Paul L\'evy in \cite{Levy1953}. Ciesielski and Taylor's formula is
the following. If $0<j_{\nu,1}<j_{\nu,2}<\ldots$ are consecutive
zeros of the Bessel function $J_{\nu}$ of the first kind, 
\[
J_{\nu}(y)=\rbr{\frac{y}{2}}^{\nu}\sum_{m=0}^{+\infty}\frac{(-1)^{m}}{m!\Gamma\rbr{m+1+\nu}}\rbr{\frac{y}{2}}^{2m},
\]
then for $t\ge0$
\begin{align}
\P\rbr{\tau^{c}>t} & =\sum_{n\ge 1}\xi_{\nu,n}e^{-j_{\nu,n}^{2}t/\rbr{2c^{2}}}\label{eq:Leevy}
\end{align}
where 
\[
\xi_{\nu,n}=\frac{\rbr{j_{\nu,n}/2}^{\nu-1}}{\Gamma\rbr{\nu+1}J_{\nu+1}\rbr{j_{\nu,n}}},\quad n=1,2,\ldots.
\]
L\'{e}vy found that the formula for $\P\rbr{\tau^{c}>t}$ has the form
(\ref{eq:Leevy}) but did not evaluate the coefficients $\xi_{\nu,k}$.
The main drawback of formula (\ref{eq:Leevy}) is that the coefficients
$\xi_{\nu,n}$ osscilate making difficult to apply it for moderate
values of $t$. In \cite{Kent:1980} Kent considered laws of hitting
times of general diffusions and found out that for broad class of
diffusions the laws of their hitting times are infinite convolutions
of elementary mixtures of exponential distributions. (Elementary mixture
of an exponential distribution is a law of a product of a $0-1$ random
variable and an independent exponential random variable). Then he
considered Bessel process and generalized Ciesielski and Taylor's
formula (not mentioning their contibution) for the case $\nu>-1$,
$x_{0}=\sqrt{z_{0}}\in[0,c)$:
\[
\P\rbr{\tau^{c}\in\dd t}= \rbr{ \frac{c^{\nu-2}}{x_{0}^{\nu}}\sum_{n\ge 1}\frac{J_{\nu}\rbr{j_{\nu,n}x_{0}/c}j_{\nu,n}}{J_{\nu+1}\rbr{j_{\nu,n}}}e^{-j_{\nu,n}^{2}t/\rbr{2c^{2}}} } \dd t.
\]
Kent's work was continued in \cite{Yamazato:1990},
where full characterisation of  distributions of hitting times (given
that they are finite) of one-dimensional diffusions (transformed so
that they are on natural scale) was given. In \cite{Yamazato:1990}
there was also the tail behaviour of these distributions characterised.
It is known that such distributions have densities which have unique
maximum (\cite{Roesler:1980}). In \cite{jedidi:2015} a characterisation
of those diffusions whose hitting times have ``bell-shaped'' distributions
was given. 

Estimates of densities of hititng times $\tau^{c}$ up to constants
depending on the order $\nu$ were given in \cite{Serafin:2017}.
In the same paper there was also established exact asyptotics of the
densities $\P\rbr{\tau^{c}\in\dd t}$ for $t\approx0$. It seems that
the literature on analytic formulas and density estimates of hitting
times $\tau^{c}$ when the starting point of the Bessel process $x_{0}$
belongs to the interval $(c,+\ns)$ is larger, see \cite{Byczkowskietal:2006,Byczkowskietal:2007,Byczkowskietal:2012}, \cite{YamanaMatsumoto:2012,YamanaMatsumoto:2013}.

In this paper we give explicit upper and lower bounds for the left ($\P \rbr{\tau^c \le t}$) 
and right ($\P \rbr{\tau^c \ge t}$)  tails of $\tau^{c}$ when the starting point $x_{0}$ of the Bessel
process belongs to the interval $[0,c)$, and exact (up to
universal, multiplicative constants, which do not dependent on $p$ neither on parameters $\nu>-1, c>0, x_0$) estimates of the ordinary moments $\Vert \tau^c\Vert_p$ and central moments $\Vert \tau^c - \E \tau^c \Vert_p$ for $p\ge 2$. Estimates of moments of even more general stopping times for Bessel processes were considered in \cite{deBlassie1987}, but in different context (their relation to moments of stopping places). Also, a very interesting invariance formula for stopping times for Bessel processes was recently established in \cite{Wisnie2018}, which can lead to some estimates of moments or more complicated functionals of stopping times of Bessel processes.

Our aim was also to provide best possible upper and lower bounds of the right and left tails of $\tau^c$. We managed to obtain upper and lower estimates of $\P \rbr{\tau^c \ge t}$ for $t \ge \gamma_0 \E \tau^c$ ($\gamma_0 \in (0, +\ns)$ is some universal constant, independent from $\nu, c, x_0$) which differ by universal multiplicative constants, that is, we found some explicit function $F$ of $t$ and parameters $\nu, c, z_0$ such that there exist some universal (independent from $\nu, c, x_0$) constants $ \gamma_0, \gamma_1, \gamma_1,  c_1, c_2 \in (0,+\ns)$ with the property that for any $t \ge \gamma_0 \E \tau^c$, 
\[
\P \rbr{\tau^c \ge  \gamma_1 t} \le c_1 \cdot F(t; \nu, c, z_0), \quad \P \rbr{\tau^c \ge  \gamma_2 t} \ge c_2 \cdot F(t; \nu, c, z_0).
\]
Universal estimates of central and ordinary moments as well as of right tails of $\tau^c$ were derived from corresponding estimates for more general variables, so called elementary mixtures of exponential distributions, which constitute broad class of diffusions hitting times, see \cite{Kent:1980}. 

Finally, we managed to obtain upper and lower estimates of $\P \rbr{\tau^c \le t}$ for $t < \gamma_0^{-1} \E \tau^c$, logarithms of which differ by universal multiplicative constants. 

Both cases are of interest especially for integer values of $\delta$, since they correspond to the exit times of a $\delta$-dimensional standard Brownian motion
from a ball with the center at $0$ and radius $c$, starting at some
point $b_{0}$ such that $\left|b_{0}\right|=x_{0}$. Thus, the obtained
estimates translate immediately into estimates for (left and right)
tails of exit times of multidimensional standard Brownian motion from
balls centered at $0$ or even more complicated regions. It is well
known that exit times of L\'{e}vy processes from bounded (open)
regions have exponentially light tails, see for example \cite[Proposition 5.2]{Grzywnyetal:2019}.
For the case of multidimensional Brownian motion see for example \cite{Peres:2010}. 
The case of left tails is important in the Large Deviations Theory, see subsection \ref{left_tail_up_b}.
Here we give bounds with explicit constants. 

Let us comment on the organisation of the paper. The second section is devoted to the estimation of right tails of $\tau^{c}$. We start with easy to obtain upper and lower bounds of the right tails (Subsections \ref{subsec:Exponential-martingale-technique_right_tails} and \ref{subsec:Kent0_left_tails}). Next, in Subsection  \ref{subsec:Kent1} we use results of Gluskin and Kwapie\'n \cite{GluskinKwapien1995} on estimates of moments of linear combinations of independent, symmetric random variables with logarithmically concave tails, to obtain estimates of moments and right tails of convolutions of elementary mixtures of exponential distributions  (Subsubsection \ref{sssect:auxiliary}). Then (Subsubsection \ref{appl_Bessel_r}) we apply estimates obtained in Subsubsection \ref{sssect:auxiliary} in the special case of the hitting times of a Bessel process. The third section is devoted to the estimation of left tails of $\tau^{c}$. We start with easy to obtain (using martingale techniques) upper (for any $\delta>0$) and lower (for $\delta >1$) bounds (Subsections \ref{left_tail_up_b} and \ref{Laplace_left_t_l_b}), and then (Subsection \ref{Kent_left_t_l_b})  we use Kent's representation to obtain upper and lower bounds of the left tails of $\tau^c$ also for $\delta \le 1$, there we also strengthen the previously obtained estimates for $\delta \approx 1$. Last, short subsection is devoted to application of the obtained estimates to right and left tails of exit times of multidimensional standard Brownian motion from balls and more complicated regions.

\section{Estimtes of right tails of hitting times $\tau^{c}$ of Bessel processes}

In this section we present upper and lower bounds for right tails
of hitting times $\tau^{c}$ . They will be obtained by application
of different techniques. First, we will apply exponential martingale
technique to obtain upper bound and representation as infinite convolution
of elementary mixtures of exponential distributions to obtain lower
bound. Next, we will present universal estimates of moments of the
hitting times (up to universal constants, independent from the dimension
of the process or the starting point) and obtain universal estimates for the right tails. 

\subsection{Exponential martingale technique and an upper bound for the right
tails}  \label{subsec:Exponential-martingale-technique_right_tails}

Let $X_{t}$, $Z_{t}$, $t\ge0$, be defined as in Section \ref{sec:Introduction}.
By the It\^{o} lemma, for any $\lambda\in\R$, the process $Y_{t}=\exp\rbr{\lambda\rbr{Z_{t}-\delta\cdot t}}$
satisfies the following stochastic differential equation:
\[
Y_{t}=Y_{0}+2\lambda\int_{0}^{t}Y_{s}\sqrt{Z_{s}}\dd{\beta_{s}}+2\lambda^{2}\int_{0}^{t}Y_{s}Z_{s}\dd s.
\]
From this we get
\begin{align*}
\E Y_{t\wedge\tau^{c}} &  =Y_{0}+2\lambda^{2}\E\int_{0}^{t\wedge\tau^{c}}Y_{s}Z_{s}\dd s  \\
& \le Y_{0}+2\lambda^{2}c^{2}\E\int_{0}^{t\wedge\tau^{c}}Y_{s\wedge\tau^{c}}\dd s 
 \le Y_{0}+2\lambda^{2}c^{2}\int_{0}^{t}\E Y_{s\wedge\tau^{c}}\dd s
\end{align*}
and by Gronwall's lemma, for any $t>0$
\[
\E Y_{t\wedge\tau^{c}}\le Y_{0}\exp\left(2\lambda^{2}c^{2}t\right)=\exp\left(2\lambda^{2}c^{2}t+\lambda\cdot z_{0}\right).
\]
This yields 
\begin{equation}
\E\exp\left(\lambda\left(Z_{t\wedge\tau^{c}}-\delta\cdot t\wedge\tau^{c}-2\lambda c^{2}t-z_{0}\right)\right)\le1.\label{eq:exp_estimate}
\end{equation}
Using (\ref{eq:exp_estimate}) for any $u\in\R$ and $\lambda>0$
we have
\begin{align*}
 & \P\rbr{\delta\cdot t-c^{2}-2\lambda c^{2}t+z_{0}=u\text{ and }\tau^{c}\ge t}\\
 & =\mathbb{P}\left(\exp\left(\lambda\left(\delta\cdot t-c^{2}-2\lambda c^{2}t+z_{0}\right)\right)\ge e^{\lambda\cdot u}\text{ and }\tau^{c}\ge t\right)\\
 & \le\mathbb{P}\left(\exp\left(\lambda\left(\delta\cdot t\wedge\tau^{c}-Z_{t\wedge\tau^{c}}-2\lambda c^{2}t+z_{0}\right)\right)\ge e^{\lambda\cdot u}\text{ and }\tau^{c}\ge t\right)\\
 & \le\E\exp\left(\lambda\left(\delta\cdot t\wedge\tau^{c}-Z_{t\wedge\tau^{c}}-2\lambda c^{2}t+z_{0}\right);\tau^{c}\ge t\right)e^{-\lambda\cdot u}\\
 & =\E\exp\left(-\lambda\left(Z_{t\wedge\tau^{c}}-\delta\cdot t\wedge\tau^{c}+2\lambda c^{2}t-z_{0}\right);\tau^{c}\ge t\right)e^{-\lambda\cdot u}\\
 & \le\E\exp\left(-\lambda\left(Z_{t\wedge\tau^{c}}-\delta\cdot t\wedge\tau^{c}-\left(-2\lambda\right)c^{2}t-z_{0}\right)\right)e^{-\lambda\cdot u} \le e^{-\lambda\cdot u}.
\end{align*}
Hence 
\begin{align*}
 \P\rbr{\delta\cdot t-c^{2}-2\lambda c^{2}t+z_{0}=u\text{ and }\tau^{c}\ge t} =\P\rbr{t=\frac{c^{2}-z_{0}+u}{\delta-2\lambda c^{2}}\text{ and }\tau^{c}\ge t} \le e^{-\lambda\cdot u}.
\end{align*}
For any $\eta>0$, substituting $u=\eta\cdot\rbr{c^{2}-z_{0}}/2$,
$\lambda=\delta\cdot\eta/(4c^{2}(1+\eta))$ we get 
\[
t=\frac{c^{2}-z_{0}+u}{\delta-2\lambda c^{2}}=\frac{c^{2}-z_{0}}{\delta}\left(1+\eta\right)
\]
and
\begin{align*}
 & \P\rbr{\tau^{c}\ge\frac{c^{2}-z_{0}}{\delta}\left(1+\eta\right)} \le\exp\rbr{-\frac{\delta\cdot\eta}{4c^{2}\left(1+\eta\right)}\frac{\eta\cdot\rbr{c^{2}-z_{0}}}{2}}=\exp\rbr{-\delta\frac{c^{2}-z_{0}}{c^{2}}\frac{\eta^{2}}{8(1+\eta)}}.
\end{align*}
Recalling that $\E\tau^{c}=\rbr{c^{2}-z_{0}}/\delta$ we can write
the just obtained estimate as 
\begin{align}
\P\rbr{\tau^{c}\ge\left(1+\eta\right)\E\tau^{c}} & \le\exp\rbr{-\delta\frac{c^{2}-z_{0}}{c^{2}}\frac{\eta^{2}}{8(1+\eta)}}.\label{eq:estim1}
\end{align}

\subsection{The law of $\tau^{c}$ as convolution of elementary mixtures of exponential
distributions and a lower bound for the right tails} \label{subsec:Kent0_left_tails}

In \cite{Kent:1980} Kent proved that the law of $\tau^{c}$ is a convolution
of elementary mixtures of exponential distributions. More precisely,
\begin{equation}
\tau^{c} \stackrel{\text{law}}{=} \sum_{n\ge 1}\tau_{\nu,n}^{c},\label{eq:law_tua_c}
\end{equation}
where $\tau_{\nu,n}^{c}$, $n=1,2,\ldots$, are independent and the
law of $\tau_{\nu,n}^{c}$ is the following: 
\begin{equation}
\P\rbr{\tau_{\nu,n}^{c}=0}=\frac{z_{0}}{c^{2}},\quad\P\rbr{\tau_{\nu,n}^{c}>t}=\rbr{1-\frac{z_{0}}{c^{2}}}\exp\rbr{-\frac{j_{\nu,n}^{2}}{2c^{2}}t},\quad t>0.\label{eq:law_tau_c1}
\end{equation}
Taking $\eta>0$ and substituting $t=\left(1+\eta\right)\E\tau^{c}=\left(1+\eta\right)\rbr{c^{2}-x_{0}^{2}}/\delta$
we have 
\begin{align}
\P\rbr{\tau^{c}\ge\left(1+\eta\right)\E\tau^{c}} & \ge\P\rbr{\tau_{\nu,1}^{c}>\left(1+\eta\right)\E\tau^{c}}=\rbr{1-\frac{z_{0}}{c^{2}}}\exp\rbr{-\frac{j_{\nu,1}^{2}}{2\delta}\frac{c^{2}-z_{0}}{c^{2}}\left(1+\eta\right)}.\label{eq:estim2}
\end{align}
Estimates (\ref{eq:estim1}) and (\ref{eq:estim2}) are of the same
order (up to a universal multiplicative constant) for large values
of $\nu$ and small values of $z_{0}$ (say $z_{0}\in\sbr{0,c^{2}/2}$).
Using for example \cite[Theorem 1 and Lemma 2]{Breen:1995} for $\nu\ge5.2$
we have 
\[
\nu+\frac{1}{7}<j_{\nu,1}<\frac{\nu}{1-3\nu^{-2/3}}
\]
and this for $\delta=2(\nu+1)\ge12.4$ yields
\begin{align}
  \P\rbr{\tau^{c}\ge\left(1+\eta\right)\E\tau^{c}}
 & \ge\rbr{1-\frac{z_{0}}{c^{2}}}\exp\rbr{-\frac{j_{\nu,1}^{2}}{2\delta}\frac{c^{2}-z_{0}}{c^{2}}\left(1+\eta\right)}\nonumber \\
 & \ge\rbr{1-\frac{z_{0}}{c^{2}}}\exp\rbr{-\frac{\nu^{2}}{2\left(1-3\nu^{-2/3}\right)^{2}\delta}\frac{c^{2}-z_{0}}{c^{2}}\left(1+\eta\right)}\nonumber \\
 & =\rbr{1-\frac{z_{0}}{c^{2}}}\exp\rbr{-\frac{\delta-4+\frac{4}{\delta}}{8\left(1-\frac{3 \cdot 2^{2/3}}{\rbr{\delta-2}^{2/3}}\right)^{2}}\frac{c^{2}-z_{0}}{c^{2}}\left(1+\eta\right)}.\label{eq:est3}
\end{align}
To summarize the estimates obtained in this and the previous subsection
we state.

\begin{proposition} \label{Right_tails_rough} Let $c>0$ and $\tau^{c}$
be the hitting time of the $\delta$-dimensional Bessel process starting
from $x_{0}=\sqrt{z_{0}}\in[0,c)$, defined by (\ref{eq:tau_c_def}).
For any $\eta>0$ the following upper and lower bounds hold: 
\[
\P\rbr{\tau^{c}\ge\left(1+\eta\right)\E\tau^{c}}\le\exp\rbr{-\delta\frac{c^{2}-z_{0}}{c^{2}}\frac{\eta^{2}}{8(1+\eta)}}
\]
and
\[
\P\rbr{\tau^{c}\ge\left(1+\eta\right)\E\tau^{c}}\ge\rbr{1-\frac{z_{0}}{c^{2}}}\exp\rbr{-\frac{j_{\nu,1}^{2}}{2\delta}\frac{c^{2}-z_{0}}{c^{2}}\left(\eta+1\right)},
\]
where $\nu=\rbr{\delta/2}-1$ and $j_{\nu,1}$ is the first zero of
the Bessel function $J_{\nu}$ of the first kind. From the last estimate
it also follows that for $\delta\ge12.4$ 
\[
\P\rbr{\tau^{c}\ge\left(1+\eta\right)\E\tau^{c}}\ge\rbr{1-\frac{z_{0}}{c^{2}}}\exp\rbr{-\frac{\delta-4+\frac{4}{\delta}}{8\left(1-\frac{3 \cdot 2^{2/3}}{\rbr{\delta-2}^{2/3}}\right)^{2}}\frac{c^{2}-z_{0}}{c^{2}}\left(\eta+1\right)}.
\]
\end{proposition}
\begin{remark}
We can naturally use the estimate $\P\rbr{\tau^{c}\ge\left(1+\eta\right)\E\tau^{c}}\ge\P\rbr{\tau_{\nu,1}^{c}+\tau_{\nu,2}^{c}\ge\left(1+\eta\right)\E\tau^{c}}$
etc. (and \cite[Theorem 1]{Breen:1995}) to obtain tighter bounds
than (\ref{eq:est3}).
Conditioning on the event that $\tau_{\nu,k}^{c} = 0$ for $k=1,2,\ldots, n-1$ and $\tau_{\nu,n}^{c} > 0$ ($n=1,2,\ldots$) we also get
\begin{align*}
  \P\rbr{\tau^{c}\ge\left(1+\eta\right)\E\tau^{c}}
 &=\sum_{n\ge1} \rbr{\frac{z_{0}}{c^{2}}}^{n-1} \rbr{1-\frac{z_{0}}{c^{2}}}\P\rbr{ \sum_{k \ge n} \tau_{\nu,k}^{c}  \ge\left(1+\eta\right)\E\tau^{c}} \nonumber \\
 & \ge \sum_{n\ge1} \rbr{\frac{z_{0}}{c^{2}}}^{n-1} \rbr{1-\frac{z_{0}}{c^{2}}}\P\rbr{ \tau_{\nu,n}^{c}  \ge\left(1+\eta\right)\E\tau^{c}}  \nonumber \\
 & =\sum_{n\ge1} \rbr{\frac{z_{0}}{c^{2}}}^{n-1} \rbr{1-\frac{z_{0}}{c^{2}}}\exp\rbr{-\frac{j_{\nu,n}^{2}}{2\delta}\frac{c^{2}-z_{0}}{c^{2}}\left(1+\eta\right)} \nonumber \\
 & \ge \exp\rbr{-\sum_{n\ge1} \rbr{\frac{z_{0}}{c^{2}}}^{n-1} \frac{j_{\nu,n}^{2}}{2\delta}\rbr{\frac{c^{2}-z_{0}}{c^{2}}}^2\left(1+\eta\right)}.
\end{align*}
Pr\'{e}kopa-Leindler inequality yields that in the case $z_0 = 0$ the tails of $\tau^c$ are log-concave.
\end{remark}

\subsection{Universal estimates of moments of $\tau^{c}$ and the right tails of $\tau^{c}$ }
\label{subsec:Kent1}

Unfortunately, since $\lim_{\nu\ra-1}j_{\nu,1}^{2}/(2\delta)=1$ (recall
that $\delta=2(\nu+1)$), see \cite{Piessens:1984}, both estimates
(\ref{eq:estim1}) and (\ref{eq:estim2}) diverge for small values
of $\delta$, even when $z_{0}=0$. This is why in this subsection
we will use different methods to provide tight (up to universal multiplicative
constants) estimates of the central moments and regular moments of $\tau^{c}$ 
\[
\left\Vert \tau^{c} - \E \tau^{c} \right\Vert _{p}=\rbr{\E \left| \tau^{c} - \E \tau^{c} \right|^{p}}^{1/p}, \quad \left\Vert \tau^{c}\right\Vert _{p}=\rbr{\E\rbr{\tau^{c}}^{p}}^{1/p}.
\]
Using them and the Paley-Zygmund inequality we will obtain tight lower
bounds for probabilities $\P\rbr{\tau^{c}\ge\gamma\left\Vert \tau^{c}\right\Vert _{p}}$
($\gamma\in(0,1)$ is a universal constant independent from $\nu$,
$c^{2}$), while the upper bound will be obtained from the Chebyshev
inequality. This method may be also applied for more general diffusions
than the Bessel process. Namely it may be applied for diffusions whose
hitting times, as those considered in \cite{Kent:1980}, have laws
which are infinite convolutions of elementary mixtures of exponential
distributions. 

To avoid technical issues related to the calculation of constants,
we introduce the following notation. For any positive quantites $a,b$
(depending on one or more parameters) we write $a\apprge b$
if there exists a universal constant $C$ such that $C\cdot a\ge b$
for all possible values of $a$ and $b$. We write $a\sim b$ and say that $a$ and $b$ are \emph{comparable} if $a\apprge b$
and $b\apprge a$. Similarly, we write $a\stackrel{p}{\apprge}b$
if $a$ and $b$ depend on the parameter $p\ge1$ and there exists
a universal constant $C$ such that $C^{p}\cdot a\ge b$ for any $p\ge1$.
We write $a\stackrel{p}{\sim}b$ if $a\stackrel{p}{\apprge}b$ and
$b\stackrel{p}{\apprge}a$. 

\subsubsection{Auxiliary lemmas} \label{sssect:auxiliary}
We start with few auxiliary results which may be of
independent interest. 
\begin{lemma}\label{Lemma_mom_general_up_b_sym} Let $\kappa_{n}$ be independent
$0$-$1$ random variables such that
\[
\P\rbr{\kappa_{n}=0}=1-\P\rbr{\kappa_{n}=1}=\alpha\in[0,1)
\]
and $\theta_{n}$ independent, exponential random variables with the
expectation $1$, $n=1,2,\ldots$. If $a_1\ge a_2\ge \ldots $ are non-negative
reals such that $\sum_{n\ge 1} a^2_{n}<+\ns$ then for $p\ge 2$ 
\begin{align*}
 \left\Vert \sum_{n \ge 1} \rbr{a_n\kappa_n\theta_n-\E a_n \kappa_n \theta_n} \right\Vert_p
 & \sim\rbr{ (1-\alpha) \sum_{n \ge 1} \alpha^{n-1}\rbr{p \cdot a_n+\sqrt{p}\sqrt{(1-\alpha)\sum_{k>n}a_k^2}}^p }^{1/p} \\
 & \sim p \rbr{1-\alpha}^{1/p} \rbr{\sum_{n\ge 1}\alpha^{n-1} a_n^{p} }^{1/p}  + \sqrt{p}\sqrt{(1-\alpha)\sum_{n \ge 1}a_n^2}.
\end{align*}
\end{lemma}
\begin{proof} 
The first fact we will use is that for $p \ge 2$
\[
\left\Vert \sum_{n\ge 1} \rbr{a_n\kappa_n\theta_n-\E a_n \kappa_n \theta_n} \right\Vert_p  \sim \left\Vert \sum_{n\ge 1} a_n\kappa_n\theta_n \varepsilon_n \right\Vert_p  
\]
where $\varepsilon_n$ are independent Rademacher variables ($\P \rbr{\varepsilon_n = \pm 1} = 1/2$), independent of $\kappa_n$ and $\theta_n$, see \cite[Remark 2]{Latala1997}. 
Next, by a result due to Gluskin and Kwapie\'n \cite{GluskinKwapien1995}, for any subset $I\subset \mathbb{N}$ we have
\begin{equation} \label{GK}
\left\Vert \sum_{n\in I} a_n\varepsilon_n \kappa_n\theta_n \right\Vert_p\sim p\Vert a(I)\Vert_{\infty}
+\sqrt{p}\Vert a(I) \Vert_2,
\end{equation} 
where $a(I)=(a_n)_{n\in I}$, $\Vert a(I)\Vert_{\infty}
 = \sup_{n\in I} \left| a_n \right|$, 
$\Vert a(I) \Vert_2 = \rbr{\sum_{n\in I} a_n^2}^{1/2}$. 
Let $\tau$ be the random time defined as the smallest $n\ge 1$ such that $\kappa_n=1$, $\tau := \min\cbr{n\ge 1: \kappa_n=1} $. By \eqref{GK} and the definition of $\tau$ we get
\[
\left\Vert \sum_{n\ge 1} a_n\varepsilon_n \kappa_n\theta_n \right \Vert_p \sim
\rbr{\E \sum_{n\ge 1}{\bf 1}_{\cbr{\tau = n}}\sbr{p \cdot a_n+\sqrt{p}\rbr{\sum_{k>n}\kappa_k a_k^2}^{1/2}}^{p}}^{1/p}.
\] 
Note that the event $\cbr{\tau=n}$ is independent of $\kappa_{n+1},\kappa_{n+2},\ldots$ and $\E {\bf 1}_{\cbr{\tau = n}} = (1-\alpha)\alpha^{n-1}$.
Therefore, to prove 
\[
 \left\Vert \sum_{n \ge 1} \rbr{a_n\kappa_n\theta_n-\E a_n \kappa_n \theta_n} \right\Vert_p
 \sim\rbr{ (1-\alpha) \sum_{n \ge 1} \alpha^{n-1}\rbr{p \cdot a_n+\sqrt{p}\sqrt{(1-\alpha)\sum_{k>n}a_k^2}}^p }^{1/p} 
 \]
it suffices to show that for $n=1,2,\ldots$
\begin{equation}\label{rafal}
\rbr{\E \sbr{p \cdot a_n+\sqrt{p}\rbr{\sum_{k>n}\kappa_k a_k^2}^{1/2}}^{p}}^{1/p}
\sim p \cdot a_n+\sqrt{p} \rbr{(1-\alpha) \sum_{k>n} a_k^2}^{1/2} .
\end{equation}

The lower bound is easy, since it is a consequence of Jensen's inequality. Namely, for $p\ge 2$
\[
\rbr{\E\sbr{\rbr{\sum_{k>n}\kappa_k a_k^2}^{1/2}}^{p}}^{2/p}\ge (1-\alpha)\sum_{k>n}a_k^2.
\]

To obtain the upper bound, we estimate  
\begin{align} \label{fiirst}
\rbr{\E\sbr{\rbr{\sum_{k>n}\kappa_k a_k^2}^{1/2}}^{p}}^{2/p}\le
 \rbr{\E{ \left|\sum_{k>n}\rbr{ \kappa_k-\E \kappa_k } a_k^2 \right|^{p/2}}}^{2/p} + (1-\alpha) \sum_{k>n}a_k^2.
\end{align}
The Bernstein inequality implies that
\[
\P\rbr{\left|\sum_{k>n}(\kappa_k-\E \kappa_k) a_k^2\right|>u}
\le 2\exp\cbr{-\frac{u^2}{2\rbr{\alpha(1-\alpha)\sum_{k>n}a_k^4}+\frac{2}{3} a_n^2u }},
\]
therefore
\begin{align} \label{seecond}
\rbr{\E\left|\sum_{k>n}(\kappa_k-\E \kappa_k) a_k^2\right|^{p/2}}^{2/p}
\apprle p \cdot a_n^2 + \sqrt{p}\rbr{\alpha(1-\alpha)\sum_{k>n}a_k^4}^{1/2}.
\end{align} 
Consequently, using the estimate (recall that $a_n \ge a_{n+1} \ge \ldots $)
\[
\sqrt{p}\rbr{\alpha(1-\alpha)\sum_{k>n}a_k^4}^{1/2} \le \sqrt{p} a_n \rbr{(1-\alpha)\sum_{k>n}a_k^2}^{1/2},
\]
which yields 
\[
p \cdot a_n^2 + \sqrt{p}\rbr{\alpha(1-\alpha)\sum_{k>n}a_k^4}^{1/2} +  (1-\alpha)\sum_{k>n}a_k^2  \apprle p \cdot a_n^2
 + (1-\alpha)\sum_{k>n}a_k^2,
\]
by \eqref{fiirst}, \eqref{seecond} we obtain 
\[
\rbr{\E\sbr{\rbr{\sum_{k>n}\kappa_k a_k^2}^{1/2}}^{p}}^{2/p} \apprle
p \cdot a_n^2
 + (1-\alpha)\sum_{k>n}a_k^2 \sim \rbr{\sqrt{p}a_n+\sqrt{(1-\alpha)\sum_{k>n}a_k^2}}^2
\]
and hence
\[
\rbr{\E \sbr{\rbr{\sum_{k>n}\kappa_k a_k^2}^{1/2}}^{p}}^{1/p}
\apprle \sqrt{p}a_n+\sqrt{(1-\alpha)\sum_{k>n}a_k^2}.
\]
It proves that
\[
\rbr{\E \sbr{pa_n+\sqrt{p}\rbr{\sum_{k>n}\kappa_k a_k^2}^{1/2}}^{p}}^{1/p}\apprle {pa_n+\sqrt{p}\sqrt{(1-\alpha)\sum_{k>n}a_k^2}},
\]
which finally establishes (\ref{rafal}). 

To finish the proof we need to prove the relation
\begin{align}
 &\rbr{ (1-\alpha) \sum_{n \ge 1} \alpha^{n-1}\rbr{p \cdot a_n+\sqrt{p}\sqrt{(1-\alpha)\sum_{k>n}a_k^2}}^p }^{1/p} \nonumber \\
 & \sim p \rbr{1-\alpha}^{1/p} \rbr{\sum_{n\ge 1}\alpha^{n-1} a_n^{p} }^{1/p}  + \sqrt{p}\sqrt{(1-\alpha)\sum_{n \ge 1}a_n^2}. \label{relll}
\end{align}
We have
\begin{align*}
 &\rbr{ (1-\alpha) \sum_{n \ge 1} \alpha^{n-1}\rbr{p \cdot a_n+\sqrt{p}\sqrt{(1-\alpha)\sum_{k>n}a_k^2}}^p }^{1/p} \\
 & \sim \rbr{ (1-\alpha) \sum_{n \ge 1} \alpha^{n-1}\cbr{\rbr{p \cdot a_n}^p + \rbr{\sqrt{p}\sqrt{(1-\alpha)\sum_{k>n}a_k^2}}^p }}^{1/p} \\
 & \apprle \rbr{ (1-\alpha) \sum_{n \ge 1} \alpha^{n-1}\cbr{\rbr{p \cdot a_n}^p + \rbr{\sqrt{p}\sqrt{(1-\alpha)\sum_{k \ge 1}a_k^2}}^p }}^{1/p} \\
 & =  \rbr{ (1-\alpha) \sum_{n \ge 1} \alpha^{n-1}\rbr{p \cdot a_n}^p + \rbr{\sqrt{p}\sqrt{(1-\alpha)\sum_{n \ge 1}a_k^2}}^p }^{1/p} \\
 & \sim p \rbr{1-\alpha}^{1/p} \rbr{\sum_{n\ge 1}\alpha^{n-1} a_n^{p} }^{1/p}  + \sqrt{p}\sqrt{(1-\alpha)\sum_{n \ge 1}a_n^2}.
 \end{align*}
On the other hand, using monotonicity of moments, we get 
\begin{align*}
 &\rbr{ (1-\alpha) \sum_{n \ge 1} \alpha^{n-1}\rbr{p \cdot a_n+\sqrt{p}\sqrt{(1-\alpha)\sum_{k>n}a_k^2}}^p }^{1/p} \\
 & \apprge  \rbr{ (1-\alpha) \sum_{n \ge 1} \alpha^{n-1}\rbr{p \cdot a_n+\sqrt{p}\sqrt{(1-\alpha)\sum_{k>n}a_k^2}}^2 }^{1/2} \\
 & \sim \rbr{ (1-\alpha) \sum_{n \ge 1} \alpha^{n-1}\rbr{p^2 \cdot a_n^2+ p (1-\alpha)\sum_{k>n}a_k^2}}^{1/2} \\
 & \apprge \rbr{ (1-\alpha)  \sum_{n \ge 1} \alpha^{n-1}p \rbr{ a_n^2+ (1-\alpha)  \sum_{k>n}a_k^2}}^{1/2}   = \sqrt{p}\sqrt{(1-\alpha)\sum_{n \ge 1}a_n^2}.
 \end{align*}
From last two estimates and the obvious relation \begin{align*}
 \rbr{ (1-\alpha) \sum_{n \ge 1} \alpha^{n-1}\rbr{p \cdot a_n+\sqrt{p}\sqrt{(1-\alpha)\sum_{k>n}a_k^2}}^p }^{1/p} \apprge  \rbr{ (1-\alpha) \sum_{n \ge 1} \alpha^{n-1}\rbr{p \cdot a_n}^p }^{1/p}.
 \end{align*}
we get  \eqref{relll}
\end{proof}
\begin{remark} \label{p2momentsym}
Direct calculation gives
\begin{align*}
 \left\Vert \sum_{n \ge 1} \rbr{a_n\kappa_n\theta_n-\E a_n \kappa_n \theta_n} \right\Vert_2 
 & = \sqrt{ (1-\alpha)(1+ \alpha) \sum_{n \ge 1} a_n^2 } \sim \sqrt{(1-\alpha)\sum_{n \ge 1}a_n^2}
\end{align*}
so from Lemma \ref{Lemma_mom_general_up_b_sym} we also infer that 
\begin{align*}
& \left\Vert \sum_{n \ge 1} \rbr{a_n\kappa_n\theta_n-\E a_n \kappa_n \theta_n} \right\Vert_p \\
 & \sim p \rbr{1-\alpha}^{1/p} \rbr{\sum_{n\ge 1}\alpha^{n-1} a_n^{p} }^{1/p}  + \sqrt{p}\left\Vert \sum_{n \ge 1} \rbr{a_n\kappa_n\theta_n-\E a_n \kappa_n \theta_n} \right\Vert_2 .
\end{align*}
\end{remark}
\begin{corollary}\label{Lemma_mom_general_up_b} Let $\kappa_{n}$ and $\theta_{n}$, $n=1,2,\ldots$ be as in Lemma 
\ref{Lemma_mom_general_up_b_sym}. If $a_1\ge a_2\ge \ldots $ are non-negative
reals such that $\sum_{n\ge 1}a_{n}<+\ns$ then for any $p\ge2$ 
\begin{align*}
\left\Vert \sum_{n\ge 1}a_{n}\kappa_{n}\theta_{n}\right\Vert _{p} & \sim p \rbr{1-\alpha}^{1/p} \rbr{\sum_{n\ge 1}\alpha^{n-1} a_n^{p} }^{1/p} + (1-\alpha) \sum_{n \ge 1} a_n.
\end{align*}
\end{corollary}
\begin{proof} Since $\sum_{n\ge 1}a_{n}<+\ns$ the expectation of $\sum_{n\ge 1}a_{n}\kappa_{n}\theta_{n}$ is finite and equals
$
\E \sum_{n\ge 1}a_{n}\kappa_{n}\theta_{n} = (1-\alpha) \sum_{n \ge 1} a_n
$. 
By Lemma \ref{Lemma_mom_general_up_b_sym} we obtain 
\begin{align*}
\left\Vert \sum_{n \ge 1} a_n\kappa_n\theta_n \right\Vert_p &    \sim \left\Vert \sum_{n \ge 1} \rbr{a_n\kappa_n\theta_n-\E a_n \kappa_n \theta_n} \right\Vert_p + \E \sum_{n \ge 1} a_n \kappa_n \theta_n\\
& \sim p \rbr{1-\alpha}^{1/p} \rbr{\sum_{n\ge 1}\alpha^{n-1} a_n^{p} }^{1/p}  + \sqrt{p}\sqrt{(1-\alpha)\sum_{n \ge 1}a_n^2} + (1-\alpha) \sum_{n \ge 1} a_n \\
& \apprge p \rbr{1-\alpha}^{1/p} \rbr{\sum_{n\ge 1}\alpha^{n-1} a_n^{p} }^{1/p}  + (1-\alpha) \sum_{n \ge 1} a_n .
\end{align*}

On the other hand, for $n=1,2,\ldots$ we have 
\[
\sqrt{p}\sqrt{(1-\alpha)\sum_{k>n}a_k^2} \apprle \sqrt{p \cdot  a_n}\sqrt{(1-\alpha)\sum_{k>n}a_k } \apprle  p \cdot a_n+(1-\alpha)\sum_{k>n}a_k
\]
which yields
\begin{align*}
   \left\Vert \sum_{n \ge 1} \rbr{a_n\kappa_n\theta_n-\E a_n \kappa_n \theta_n} \right\Vert_p
& \sim\rbr{ (1-\alpha) \sum_{n \ge 1} \alpha^{n-1}\rbr{p \cdot a_n+\sqrt{p}\sqrt{(1-\alpha)\sum_{k>n}a_k^2}}^p }^{1/p} \\
& \apprle \rbr{ (1-\alpha) \sum_{n \ge 1} \alpha^{n-1}\rbr{p \cdot a_n+(1-\alpha) \sum_{k >n } a_n}^p }^{1/p}.
\end{align*}
To conclude we need to prove the relation 
\begin{align*}
\rbr{\rbr{1-\alpha}\sum_{n\ge 1}\alpha^{n-1}\rbr{p\cdot a_n + (1-\alpha)\sum_{k>n}a_{k}}^{p} }^{1/p} \\
\apprle p \rbr{1-\alpha}^{1/p} \rbr{\sum_{n\ge 1}\alpha^{n-1} a_n^{p} }^{1/p} + (1-\alpha) \sum_{n \ge 1} a_n.
\end{align*}
but it may be proved in a similar way as relation \eqref{relll}.
\end{proof}
Easy consequences of Corollary \ref{Lemma_mom_general_up_b} are two following corollaries.
\begin{corollary}\label{Corollary_mom_sigma} Under the assumptions of Corollary 
\ref{Lemma_mom_general_up_b} for any $p\ge2$ the $p$th moment of the variable 
\[
\sigma=a_{1}\theta_{1}+\sum_{n >1}a_{n}\kappa_{n}\theta_{n}
\]
may be estimated as 
\[
\left\Vert \sigma\right\Vert _{p}\sim  p \cdot a_{1}+(1-\alpha)\sum_{n>1}a_{n}.
\]
\end{corollary}
\begin{proof} Using Corollary \ref{Lemma_mom_general_up_b}, equality $\left\Vert \theta_{1}\right\Vert _{p} = \rbr{\Gamma(p+1)}^{1/p} \sim p $ and the estimate  
\[
\rbr{1-\alpha}^{1/p} \rbr{\sum_{n >1}\alpha^{n-2} a_n^{p} }^{1/p} \apprle   \rbr{1-\alpha}^{1/p} \rbr{\sum_{n >1}\alpha^{n-2} a_1^{p} }^{1/p} = a_1   
\]
we have 
\begin{align*}
\left\Vert \sigma\right\Vert _{p} & \sim  \left\Vert a_{1}\theta_{1}\right\Vert _{p}+\left\Vert \sum_{n > 1} a_{n}\kappa_{n}\theta_{n}\right\Vert _{p}  \\ & \sim p \cdot a_{1} +p \rbr{1-\alpha}^{1/p} \rbr{\sum_{n > 1}\alpha^{n-2} a_n^{p} }^{1/p} + (1-\alpha) \sum_{n > 1} a_n \\
 & \sim p \cdot a_{1} + (1-\alpha)\sum_{n > 1} a_n.
\end{align*}
\end{proof}
\begin{corollary}\label{Corollary_tail_sigma} Under the assumptions of Corollary
\ref{Corollary_mom_sigma} and additional assumption $a_1 >0$ there exist universal constants $ \beta \in [1,+\ns)$ and $\gamma_{10}, \gamma_{11} \in (0, +\ns)$, such that for $p\ge 2$ we have the following tail estimates: 
\begin{equation}
\P\rbr{\sigma\ge  \gamma_{10} \rbr{p \cdot a_{1} + (1-\alpha)\sum_{n > 1} a_n }} \le e^{-p}\label{eq:tail_sigma_above}
\end{equation}
and
\begin{equation}
\P\rbr{\sigma \ge \gamma_{11} \rbr{p \cdot a_{1} + (1-\alpha)\sum_{n > 1} a_n }  }\ge \frac{1}{2}e^{-\beta \cdot p}.\label{eq:sigma_tail_below}
\end{equation}
\end{corollary}
\begin{proof}
By Corollary \ref{Corollary_mom_sigma} there exist universal constants $\gamma_{12}, \gamma_{13}$ such that for $p\ge 2$, 
\begin{equation} \label{corrr213}
\gamma_{12} \rbr{p \cdot a_{1} + (1-\alpha)\sum_{n > 1} a_n } \le \left\Vert \sigma\right\Vert _{p} \le \gamma_{13} \rbr{p \cdot a_{1} + (1-\alpha)\sum_{n > 1} a_n }.
\end{equation}
This and Chebyschev's inequality yield (\ref{eq:tail_sigma_above}) with $\gamma_{10} = e \cdot \gamma_{13}$:
\begin{align*}
& \P\rbr{\sigma \ge e\cdot \gamma_{13} \rbr{p \cdot a_{1} + (1-\alpha)\sum_{n >1 } a_n }} \le \P\rbr{\sigma \ge e\left\Vert \sigma\right\Vert _{p}  } \le \frac{\E \sigma^p}{ e^p \left\Vert \sigma\right\Vert _{p}^p}  =e^{-p}.
\end{align*}
Inequality (\ref{eq:sigma_tail_below}) follows from \eqref{corrr213}
and the Paley-Zygmund inequality: for any $p\ge2$ and $\beta_0\in[0,1]$
\begin{align*}
&  \P\rbr{\sigma\ge \beta_0 \cdot \gamma_{12} \rbr{p \cdot a_{1} + (1-\alpha)\sum_{n > 1} a_n } } \\
& \ge \P\rbr{\sigma\ge\beta_0\left\Vert \sigma\right\Vert _{p}} 
 =\P\rbr{\sigma^{p}\ge\beta_0^{p}\E \sigma^p}\ge\rbr{1-\beta_0^{p}}^{2}\frac{\left\Vert \sigma\right\Vert _{p}^{2p}}{\left\Vert \sigma\right\Vert _{2p}^{2p}}\\
 & \ge\rbr{1-\beta_0^{p}}^{2}\frac{\rbr{\gamma_{12} \left\{ p\cdot a_{1}+(1-\alpha)\sum_{n=2}^{+\ns}a_{n}\right\} }^{2p}}{\rbr{\gamma_{13}\left\{ 2p \cdot a_{1}+(1-\alpha)\sum_{n=2}^{+\ns}a_{n}\right\} }^{2p}}\\
 & \ge\rbr{1-\beta_0^{p}}^{2}\rbr{\frac{\gamma_{12}}{2 \gamma_{13}}}^{2p} = \rbr{1-\beta_0^{2}}^{2}\exp\rbr{- 2 \rbr{ \ln\frac{2\gamma_{13}}{ \gamma_{12}}} p }.
\end{align*}
Substituting $\beta_0=\rbr{1 - 2^{-1/2}}^{1/2}$  we get (\ref{eq:sigma_tail_below}) with $\gamma_{11} = \beta_0 \cdot \gamma_{12}$ and $\beta = 2 \rbr{ \ln\frac{2\gamma_{13}}{ \gamma_{12}}} \ge 2 \ln 2 \ge 1$.
\end{proof}
From the last corollary we will obtain upper and lower estimates of
right tails of the variables $\sum_{n\ge 1}a_{n}\kappa_{n}\theta_{n}$. 
\begin{proposition}\label{Proposition.tails_tau_gen} For any $t > 0$ and
non-negative $a_{1}\ge a_{2}\ge a_{3}\ge\ldots$ such that  $\sum_{n\ge 1} a_n <+\ns$ let $p_{n}(t)$, $n=1,2,\ldots$, be defined by the equality 
\[
t=  p_{n}(t)\cdot a_{n}+(1-\alpha)\sum_{k>n}a_{k} 
\]
 if $a_n >0$ and $p_{n}(t) = + \ns$  if $a_n = 0$.
If $\kappa_{n}, \theta_{n}$, $n=1,2,\ldots,$ are as in Lemma \ref{Lemma_mom_general_up_b_sym}  then there exist universal constants $\beta \ge 1$ and $\gamma_7, \gamma_8, \gamma_9 \in (0,+\ns)$ such that we have the following estimates 
\begin{equation} \label{og_gorne}
\P\rbr{\sum_{n\ge 1}a_{n}\kappa_{n}\theta_{n}\ge \gamma_7 \cdot t}\le (1-\alpha) \cbr{\sum_{n:p_{n}(t)<2}\alpha^{n-1}+\sum_{n:p_{n}(t)\ge 2}\alpha^{n-1}e^{-\beta \cdot p_{n}(t)}}
\end{equation}
and 
\begin{equation} \label{og_dolne}
\P\rbr{\sum_{n\ge 1}a_{n}\kappa_{n}\theta_{n}\ge \gamma_8 \cdot t}\ge 
\gamma_9 (1-\alpha) \cbr{  \sum_{n:p_{n}(t)<2}\alpha^{n-1}+ \sum_{n:p_{n}(t)\ge2}\alpha^{n-1} e^{-\beta \cdot p_{n}(t)} }.
\end{equation}
\end{proposition}
\begin{proof} 
Let $\tau$ be the random time defined as in the proof of Lemma \ref{Lemma_mom_general_up_b_sym}, that is as the smallest $n\ge 1$ such that $\kappa_n=1$, $\tau := \min\cbr{n\ge 1: \kappa_n=1} $ and let $\sigma_{n}$, $n=1,2,\ldots$, denote the variable 
\[
\sigma_{n}=a_{n}\theta_{n}+\sum_{k > n} a_{k}\kappa_{k}\theta_{k}.
\]
By Corollary \ref{Corollary_tail_sigma},  estimate \eqref{eq:sigma_tail_below}, 
we have 
\[
\P\rbr{\sigma_{n}\ge \gamma_{11} t}\ge\begin{cases}
\frac{1}{2}e^{- \beta \cdot 2} & \text{if }p_{n}(t)<2,\\
\frac{1}{2}e^{-\beta \cdot p_{n}(t)} & \text{if }p_{n}(t)\ge 2,
\end{cases}
\]
for some universal $\beta \ge1$, from which we estimate 
\begin{align*}
  \P\rbr{\sum_{n\ge 1}a_{n}\kappa_{n}\theta_{n} \ge \gamma_{11} t}& =\sum_{n\ge 1} \P\rbr{\sum_{k \ge 1} a_{k}\kappa_{k}\theta_{k} \ge \gamma_{11} t | \tau = n}\P\rbr{\tau = n}\\
 & =\sum_{n\ge 1}\P\rbr{\sigma_{n} \ge \gamma_{11} t}\P\rbr{\tau = n}\\
 & \ge\sum_{n:p_{n}(t)<2}\frac{1}{2}e^{- \beta \cdot 2}(1-\alpha)\alpha^{n-1}+\sum_{n:p_{n}(t)\ge2}\frac{1}{2}e^{- \beta \cdot p_{n}(t)}(1-\alpha)\alpha^{n-1},
\end{align*}
which gives \eqref{og_dolne} with $\gamma_8 = \gamma_{11}$, $\gamma_9 =\frac{1}{2}e^{- \beta \cdot 2}$.
On the other hand, if $a_n >0$ we have 
\[
p_{n}(\beta t) = \frac{1}{a_n} \beta t - \frac{1}{a_n} (1-\alpha)\sum_{k > n}a_{k} \ge \beta \rbr{\frac{1}{a_n} t - \frac{1}{a_n} (1-\alpha)\sum_{k > n}a_{k}} = \beta \cdot p_{n}(t).
\]
which together with \eqref{eq:tail_sigma_above} gives
\begin{align*}
\P\rbr{\sigma_{n}\ge \gamma_{10}\beta t} & =\P\rbr{\sigma_{n}\ge \gamma_{10} \rbr{ p_{n}(\beta t)\cdot a_{n}+(1-\alpha)\sum_{k>n}a_{k}}}\\
 & \le \P\rbr{\sigma_{n}\ge \gamma_{10} \rbr{ \beta \cdot p_{n}(t)\cdot a_{n}+(1-\alpha)\sum_{k>n}a_{k}}}\\ 
 & \le\begin{cases}
1 & \text{ if } \beta \cdot p_{n}(t) <2,\\
e^{-\beta \cdot p_{n}(t)} & \text{ if }\beta \cdot p_{n}(t)\ge2 
\end{cases}
\le\begin{cases}
1 & \text{ if } p_{n}(t) <2,\\
e^{-\beta \cdot p_{n}(t)} & \text{ if } p_{n}(t)\ge2 
\end{cases}
\end{align*}
and we estimate 
\begin{align*}
  \P\rbr{\sum_{n\ge 1}a_{n}\kappa_{n}\theta_{n} \ge \gamma_{10} \beta t}& =\sum_{n\ge 1} \P\rbr{\sum_{k \ge 1}a_{k}\kappa_{k}\theta_{k} \ge \gamma_{10} \beta t | \tau = n}\P\rbr{\tau = n}\\
 & =\sum_{n\ge 1}\P\rbr{\sigma_{n} \ge \gamma_{10} \beta t}\P\rbr{\tau = n}\\
 & \le\sum_{n:p_{n}(t)<2} 1\cdot (1-\alpha)\alpha^{n-1}+\sum_{n:p_{n}(t)\ge2}e^{- \beta \cdot p_{n}(t)}(1-\alpha)\alpha^{n-1},
\end{align*}
which gives \eqref{og_gorne} with $\gamma_7 = \gamma_{10} \beta$.
\end{proof}
\subsubsection{Application to Bessel processes} \label{appl_Bessel_r}
Now we are ready to obtain universal (up to universal multiplicative
constants) estimates of $\left\Vert \tau^{c} - \E \tau^c \right\Vert _{p}$, $\left\Vert \tau^{c}\right\Vert _{p}$ and of tails of $\tau^{c}$ (we apply the notation from Section \ref{sec:Introduction} and Subsection \ref{subsec:Kent0_left_tails}). 

The first theorem we prove establishes universal estimates of the central moments of $\tau^c$.
\begin{theorem} \label{MomentsCentral} Let $p\ge2$ and let us denote $\alpha=z_{0}/c^{2} = x_0^2/c^2$. The following estimates hold
\begin{itemize}
\item if $\alpha^{\nu+2}\ge1/2$ then 
\begin{equation}
\left\Vert \frac{\tau^c - \E \tau^c}{c^2} \right\Vert_p \sim p\rbr{1-\alpha}^{1/p}\cbr{\frac{1}{j_{\nu,1}^2}+ \frac{1}{(\nu+2)^{2-1/p}}}+\frac{\sqrt{p}\sqrt{1-\alpha}}{(\nu+1)\sqrt{\nu+2}} ;\label{eq:first_case_sym}
\end{equation}
\item if  $\alpha^{\nu+2}<1/2$ but $\alpha > 1/2$ then 
\begin{equation}
\left\Vert \frac{\tau^c - \E \tau^c}{c^2} \right\Vert_p \sim p\rbr{1-\alpha}^{1/p}\cbr{\frac{1}{j_{\nu,1}^2}+ \frac{1}{(\nu+2)^{2} \rbr{\log_2(1/\alpha)}^{1/p}}}+\frac{\sqrt{p}\sqrt{1-\alpha}}{(\nu+1)\sqrt{\nu+2}} \label{eq:second_case_sym}
\end{equation}
\item and if $\alpha\le1/2$ then 
\begin{equation}
\left\Vert \frac{\tau^c - \E \tau^c}{c^2} \right\Vert_p \sim p\rbr{1-\alpha}^{1/p}\frac{1}{j_{\nu,1}^2}+\frac{\sqrt{p}\sqrt{1-\alpha}}{(\nu+1)\sqrt{\nu+2}} \sim \frac{p}{j_{\nu,1}^2}+\frac{\sqrt{p}}{(\nu+1)\sqrt{\nu+2}}.\label{eq:third_case_sym}
\end{equation}
\end{itemize}
\end{theorem}
\begin{proof} 
Using \eqref{eq:law_tua_c}, \eqref{eq:law_tau_c1} and Lemma \ref{Lemma_mom_general_up_b_sym} we infer that
\begin{align*}
\left\Vert \tau^c - \E \tau^c \right\Vert_p \sim p \rbr{1-\alpha}^{1/p} \rbr{\sum_{n\ge 1}\alpha^{n-1} a_n^{p} }^{1/p}  + \sqrt{p}\sqrt{(1-\alpha) \sum_{n\ge 1} a_n^2}, 
\end{align*}
where 
$
a_n={2c^2}/{ j_{\nu,n}^{2}}.
$
Thus 
\begin{align}
\left\Vert \frac{\tau^c - \E \tau^c}{c^2} \right\Vert_p \sim p \rbr{1-\alpha}^{1/p} \rbr{\sum_{n\ge 1}\alpha^{n-1} \frac{1}{j_{\nu,n}^{2p}} }^{1/p}  + \sqrt{p}\sqrt{(1-\alpha) \sum_{n\ge 1} \frac{1}{j_{\nu,n}^{4}}}, \label{first_step}
\end{align}
To calculate $ \sum_{n\ge 1}  1/j_{\nu,n}^{4} $ we may use the identities
\begin{align} 
{\Gamma(\nu+1)} \rbr{\frac{y}{2}}^{- \nu} J_{\nu}(y) & = \prod_{n\ge 1}\rbr{1- \frac{y^2}{j_{\nu,n}^{2}}} =\sum_{m\ge 0}\frac{(-1)^m{\Gamma(\nu+1)} }{m!\Gamma\rbr{\nu+1+m}}\rbr{\frac{y}{2}}^{2m}, \label{Jexp} \\
{\Gamma(\nu+1)}  \rbr{\frac{y}{2}}^{- \nu} I_{\nu}(y) & = \prod_{n \ge 1}\rbr{1+\frac{y^2}{j_{\nu,n}^{2}}} =\sum_{m \ge 0}\frac{{\Gamma(\nu+1)} }{m!\Gamma\rbr{\nu+1+m}}\rbr{\frac{y}{2}}^{2m} \nonumber
\end{align}
($I_{\nu}$ is the modified Bessel function of the first kind), and we have
\[
 \prod_{n \ge 1}\rbr{1- \frac{y^2}{j_{\nu,n}^{2}}} \prod_{n \ge 1}\rbr{1+ \frac{y^2}{j_{\nu,n}^{2}}} =1 - \rbr{\sum_{n\ge 1}  \frac{1}{j_{\nu,n}^{4}}}y^4 + \ldots = 1 - \frac{y^4}{16(\nu+1)^2(\nu+2)} + \ldots
\]
thus 
\begin{equation} \label{first_ingr}
\sum_{n\ge 1}  \frac{1}{j_{\nu,n}^{4}}  = \frac{1}{16(\nu+1)^2(\nu+2)}.
\end{equation}

Now we will estimate the sum $\sum_{n\ge 2}\alpha^{n-1}/{j_{\nu,n}^{2p}}$. Since the zeros of Bessel functions of the first kind are interlacing,
i.e. $j_{\nu,1}<j_{\nu+1,1}<j_{\nu,2}<j_{\nu+1,2}<\ldots$ (see \cite[Chap. XV]{Watson:1944})
we notice that for $n\ge2$, $j_{\nu,n}>j_{\nu+1,n-1}$. Now, using
\cite[p. 490]{Watson:1944}) for $\nu\in(-1,-1/2)$ and \cite[Theorem 1 and Lemma 2]{Breen:1995}
for $\nu\ge-1/2$ we get that for any $\nu>-1$ and $n\ge2$ 
\begin{equation}
j_{\nu,n}\sim j_{\nu+1,n-1}\sim\nu+n.\label{eq:zeros1}
\end{equation}

We have three cases. If $\alpha^{\nu+2}\ge1/2$ then using (\ref{eq:zeros1})
we get
\begin{align}
  \sum_{n \ge 2}\alpha^{n-1}{\frac{1}{j_{\nu,n}^{2p}}} & \stackrel{p}{\sim} \sum_{n \ge 2}\alpha^{n-1}\frac{1}{\rbr{\nu+n}^{2p}}\nonumber \\
 & \ge  \sum_{n=2}^{\left\lfloor \nu\right\rfloor +3}\alpha^{n-1}\frac{1}{\rbr{\nu+n}^{2p}}\stackrel{p}{\sim}\frac{1}{2}\frac{\left\lfloor \nu\right\rfloor +2}{\rbr{\nu+2+\left\lfloor \nu\right\rfloor +2}^{2p}}\nonumber \\
 & \stackrel{p}{\sim}  \rbr{\nu+2}^{1-2p}.\label{eq:mom_jeden}
\end{align}
On the other hand,
\begin{align}
 \sum_{n\ge 2}\alpha^{n-1}\frac{1}{j_{\nu,n}^{2p}} & \stackrel{p}{\sim} \sum_{n\ge 2}\alpha^{n-1}\frac{1}{\rbr{\nu+n}^{2p}}
 \le\sum_{n \ge 2}\frac{1}{\rbr{\nu+n}^{2p}} \nonumber \\
 & \stackrel{p}{\apprle}  \int_{\nu+2}^{+\ns}\frac{1}{x^{2p}}\dd x  = \frac{1}{2p-1} \rbr{\nu+2}^{1-2p}\stackrel{p}{\sim} \rbr{\nu+2}^{1-2p}.\label{eq:mom_dwa}
\end{align}
From (\ref{eq:mom_jeden}) and (\ref{eq:mom_dwa}) we get 
\begin{equation} \label{sec_ingr}
\sum_{n \ge 2}\alpha^{n-1}\frac{1}{j_{\nu,n}^{2p}}\stackrel{p}{\sim}\rbr{\nu+2}^{1-2p}.
\end{equation}
From \eqref{first_step}, \eqref{sec_ingr} and \eqref{first_ingr}
we obtain (\ref{eq:first_case_sym}). 

Now let us consider the case $\alpha^{\nu+2}<1/2$ and $\alpha>1/2$ (that is $ 1 < 1/\log_2(1/\alpha)< \nu +2 $).
Using the relations $\alpha^{1/\log_{2}(1/\alpha)} = 1/2$ and $ 1 < 1/\log_2(1/\alpha)< \nu +2 $ we estimate 
\begin{align}
 \sum_{n \ge 2}\alpha^{n-1}\frac{1}{j_{\nu,n}^{2p}} &\stackrel{p}{\sim} \sum_{n \ge 2}\alpha^{n-1}\frac{1}{\rbr{\nu+n}^{2p}}   \ge \sum_{n=2}^{\left\lfloor 1/\log_{2}(1/\alpha)\right\rfloor +1}\alpha^{n-1}\frac{1}{\rbr{\nu+n}^{2p}}\nonumber \\
 & \ge\frac{1}{2}\frac{\left\lfloor 1/\log_{2}(1/\alpha)\right\rfloor }{\rbr{\nu+\left\lfloor 1/\log_{2}(1/\alpha)\right\rfloor +1}^{2p}} \stackrel{p}{\apprge} \frac{\rbr{\nu+2}^{-2p}}{\log_{2}(1/\alpha)}.\label{eq:mom_3}
\end{align}
On the other hand,
\begin{align}
\sum_{n \ge 2}\alpha^{n-1}\frac{1}{j_{\nu,n}^{2p}} &\stackrel{p}{\sim} \sum_{n \ge 2}\alpha^{n-1}\frac{1}{\rbr{\nu+n}^{2p}}   \nonumber \\
 & \le \sum_{k=0}^{+\ns}\sum_{n=k\left\lfloor 1+1/\log_{2}(1/\alpha)\right\rfloor +2}^{(k+1)\left\lfloor 1+1/\log_{2}(1/\alpha)\right\rfloor +1}\alpha^{k\left\lfloor 1+1/\log_{2}(1/\alpha)\right\rfloor }\frac{1}{\rbr{\nu+n}^{2p}}\nonumber \\
 & \le\sum_{k=0}^{+\ns}2^{-k}\frac{\left\lfloor 1+1/\log_{2}(1/\alpha)\right\rfloor }{\rbr{\nu+2}^{2p}} 
\apprle \frac{\rbr{\nu+2}^{-2p}}{\log_{2}(1/\alpha)}.\label{eq:mom4}
\end{align}
From \eqref{eq:mom_3} and \eqref{eq:mom4} we get 
\begin{equation} \label{mom5}
\sum_{n \ge 2}\alpha^{n-1}\frac{1}{j_{\nu,n}^{2p}} \stackrel{p}{\sim}\frac{\rbr{\nu+2}^{-2p}}{\log_{2}(1/\alpha)}.
\end{equation}
From \eqref{first_step},\eqref{mom5} and \eqref{first_ingr} we obtain \eqref{eq:second_case_sym}.

Finally, if $\alpha\le1/2$ then $\sum_{n \ge 2}\alpha^{n-1}/{j_{\nu,n}^{2p}} \apprle 1/{j_{\nu,1}^{2p}}$ and \eqref{eq:third_case_sym} follows from  \eqref{first_step} and \eqref{first_ingr}.
\end{proof}
Now we will establish universal estimates of the ordinary moments of $\tau^c$.
\begin{theorem} \label{Moments} Let $p\ge2$ and let us denote $\alpha=z_{0}/c^{2} = x_0^2/c^2$. The following estimates hold
\begin{itemize}
\item if $\alpha^{\nu+2}\ge1/2$ then 
\begin{equation}
\left\Vert \frac{\tau^{c}}{c^{2}}\right\Vert _{p}\sim p\rbr{1-\alpha}^{1/p}\cbr{\frac{1}{j_{\nu,1}^{2}}+\frac{1}{\rbr{\nu+2}^{2-1/p}}}+\frac{{1-\alpha}}{\nu+1};\label{eq:first_case}
\end{equation}
\item if $\alpha^{\nu+2}<1/2$ but $\alpha>1/2$ then 
\begin{equation}
\left\Vert \frac{\tau^{c}}{c^{2}}\right\Vert _{p}\sim p\rbr{1-\alpha}^{1/p}\cbr{\frac{1}{j_{\nu,1}^{2}}+\frac{1}{\rbr{\nu+1}^{2}\rbr{\log_2(1/\alpha)}^{1/p}}}+\frac{{1-\alpha}}{\nu+1} \label{eq:second_case}
\end{equation}
\item and if $\alpha\le1/2$ then 
\begin{equation}
\left\Vert \frac{\tau^{c}}{c^{2}}\right\Vert _{p}\sim p\rbr{1-\alpha}^{1/p}\frac{1}{j_{\nu,1}^{2}}+\frac{{1-\alpha}}{\nu+1} \sim\frac{p}{j_{\nu,1}^{2}}+\frac{1}{\nu+1}.\label{eq:third_case}
\end{equation}
\end{itemize}
\end{theorem}
\begin{proof} By (\ref{eq:law_tua_c}), (\ref{eq:law_tau_c1}) and Corollary \ref{Lemma_mom_general_up_b},
we get 
\begin{align*}
\left\Vert \tau^c \right\Vert _{p} & \sim p \rbr{1-\alpha}^{1/p} \rbr{\sum_{n\ge 1}\alpha^{n-1} a_n^{p} }^{1/p} + (1-\alpha) \sum_{n \ge 1} a_n
\end{align*}where 
$
a_n={2c^2}/{ j_{\nu,n}^{2}}.
$
Thus 
\begin{align}
\left\Vert \frac{\tau^c}{c^2} \right\Vert_p \sim p \rbr{1-\alpha}^{1/p} \rbr{\sum_{n\ge 1}\alpha^{n-1} \frac{1}{j_{\nu,n}^{2p}} }^{1/p}  + (1-\alpha) \sum_{n \ge 1} \frac{1}{j_{\nu,n}^{2}}. \label{ffirst_step}
\end{align}
The first term on the right side of \eqref{ffirst_step} was estimated in the proof of Theorem \ref{MomentsCentral}
 and since $\sum_{n\ge 1}1/{j_{\nu,n}^{2}}=1/(2\delta) = 1/(4(\nu +1))$ (this may be proved for example from \eqref{Jexp}) the result follows. 
\end{proof}
From Proposition \ref{Proposition.tails_tau_gen} we will infer the following
upper and lower estimates of right tails of $\tau^{c}$. Our method
will work for the case $t\gtrsim\E\tau^{c}$.
\begin{theorem}\label{Tails} Let us denote $\alpha=z_{0}/c^{2}=x_{0}^{2}/c^{2}$. There exist
universal constants $\gamma_{1},\gamma_{2},\gamma_{3},\gamma_{4} \in (0,+\ns)$
such that for $t\ge\gamma_{3}\E\tau^{c}$
\[
\P\rbr{\tau^{c}\ge\gamma_{1}t}\apprle (1-\alpha)\cbr{\exp\rbr{-\frac{j_{\nu,1}^{2}}{2c^{2}}t}+\sum_{n \ge 2}\alpha^{n-1}\exp\rbr{-\frac{(\nu+n)^{2}}{2c^{2}}t}}
\]
and for $t\ge\gamma_{4}\E\tau^{c}$
\[
\P\rbr{\tau^{c}\ge\gamma_{2}t} \apprge (1-\alpha)\cbr{\exp\rbr{-\frac{j_{\nu,1}^{2}}{2c^{2}}t}+\sum_{n \ge 2}\alpha^{n-1}\exp\rbr{-\frac{(\nu+n)^{2}}{2c^{2}}t}}.
\]
\end{theorem}
\begin{proof} Let us denote $a_{n}=2c^{2}/j_{\nu,n}^{2}$ for $n=1,2,\ldots$.  For $t\ge2(1-\alpha)\sum_{k=1}^{+\ns}a_{k}=2\E\tau^{c}$ let
$u = u(t)$ be defined as 
\[
u(t)=\frac{t}{2\left(1-\alpha\right)c^{2}}
\]
and $p_{n}(t)$ be defined as in Proposition \ref{Proposition.tails_tau_gen}.
By the definition of $p_{n}(t)$ and $u$ 
\begin{align}
p_{n}(t) & =\frac{t}{a_{n}}-\frac{1}{a_{n}}(1-\alpha)\sum_{k > n}a_{k} \nonumber =\frac{1}{a_{n}}\left\{ 2\left(1-\alpha\right)c^{2}u(t)-(1-\alpha)\sum_{k > n}a_{k}\right\} \nonumber \\
 & =\frac{j_{\nu,n}^{2}}{2c^{2}}(1-\alpha)\cbr{ 2 c^{2}u(t)-\sum_{k > n}\frac{2c^{2}}{j_{\nu,k}^{2}}}  =j_{\nu,n}^{2}\left(1-\alpha\right)\left\{ u(t)-\sum_{k > n}\frac{1}{j_{\nu,k}^{2}}\right\} .\label{eq:p_n_estim}
\end{align}
Since $t\ge2(1-\alpha)\sum_{k=1}^{+\ns}a_{k}$ we get
\begin{align*}
u(t) & =\frac{t}{2\left(1-\alpha\right)c^{2}}\ge\frac{1}{c^{2}}\sum_{k \ge 1}a_{k}=2\sum_{k \ge 1}\frac{1}{j_{\nu,k}^{2}},
\end{align*}
which for $n=1,2\ldots$ yields $u(t)-\sum_{k > n}1/j_{\nu,k}^{2}\sim u(t)$
and by (\ref{eq:p_n_estim}) $p_{n}(t)\sim \left(1-\alpha\right)  j_{\nu,n}^{2}\cdot u(t)$.
Recall also (see the proof of Theorem \ref{MomentsCentral}) that for $n\ge2$,
$j_{\nu,n}\sim\nu+n$, thus, for $n\ge2$, $p_{n}(t)\sim(\nu+n)^{2}\left(1-\alpha\right) \cdot u(t)$.
Let $0<\gamma_{5}\le\gamma_{6}<+\infty$ be universal constants such that 
\[
\gamma_{5}\left(1-\alpha\right) j_{\nu,1}^{2}\cdot u(t) \le p_{1}(t)\le\gamma_{6} \left(1-\alpha\right) j_{\nu,1}^{2} \cdot u(t)
\]
and such that for $n\ge2$
\[
\gamma_{5}\left(1-\alpha\right)(\nu+n)^{2}\cdot u(t)\le p_{n}(t)\le\gamma_{6}\left(1-\alpha\right) (\nu+n)^{2} \cdot u(t).
\]
Let us also notice that since $p_{n}(t)> 0$ when $t\ge2\E\tau^{2}$,
for any $n=1,2,\ldots$ such that $p_{n}(t)<2$ we have 
\[
e^{-\beta \cdot p_{n}(t)}\sim1.
\]
By this observation and Proposition \ref{Proposition.tails_tau_gen} for $t \ge 2\E\tau^{2}$
we have
\begin{align*}
 &  \P\rbr{\tau^{c}\ge\gamma_7 \cdot t} \apprle (1-\alpha)\sum_{n=1}^{+\infty}\alpha^{n-1} e^{-\beta \cdot p_{n}(t)}\\
 & \le (1-\alpha)\cbr{ \exp\rbr{-\beta  \gamma_5  (1-\alpha) j_{\nu,1}^{2}   \cdot u(t) }+\sum_{n \ge 2} \alpha^{n-1} \exp\rbr{-\beta  \gamma_5  (1-\alpha) (\nu+2)^{2}   \cdot u(t) } }.
\end{align*}
Thus, for $t\ge2 \beta \gamma_{5}\E\tau^{c}$ we have
\begin{align*}
&  \P\rbr{\tau^{c}\ge \gamma_7 \frac{t}{\beta \gamma_{5}}} \\ 
& \apprle  (1-\alpha)\cbr{ \exp\rbr{-\beta  \gamma_5  (1-\alpha) j_{\nu,1}^{2}   \cdot u\rbr{\frac{t}{\beta \gamma_{5}}} }+\sum_{n \ge 2} \alpha^{n-1} \exp\rbr{-\beta  \gamma_5  (1-\alpha) (\nu+2)^{2}   \cdot u\rbr{\frac{t}{\beta \gamma_{5}}} } } \\
 & =(1-\alpha)\cbr{\exp\rbr{-\frac{j_{\nu,1}^{2}}{2c^{2}}t}+\sum_{n \ge 2}\alpha^{n-1}\exp\rbr{-\frac{(\nu+n)^{2}}{2c^{2}}t}}.
\end{align*}
Similarly, for $t\ge2\beta \gamma_{6}\E\tau^{c}$ we have 
\begin{align*}
&  \P\rbr{\tau^{c}\ge \gamma_8 \frac{t}{\beta \gamma_{6}}} \\ 
& \apprge  (1-\alpha)\cbr{ \exp\rbr{-\beta  \gamma_6  (1-\alpha) j_{\nu,1}^{2}   \cdot u\rbr{\frac{t}{\beta \gamma_{6}}} }+\sum_{n \ge 2} \alpha^{n-1} \exp\rbr{-\beta  \gamma_6  (1-\alpha) (\nu+2)^{2}   \cdot u\rbr{\frac{t}{\beta \gamma_{6}}} } } \\
 & = (1-\alpha)\cbr{\exp\rbr{-\frac{j_{\nu,1}^{2}}{2c^{2}}t}+\sum_{n \ge 2}\alpha^{n-1}\exp\rbr{-\frac{(\nu+n)^{2}}{2c^{2}}t}}.
\end{align*}
Denoting $\gamma_{1}=\gamma_{7}/\rbr{\beta \gamma_{5}}$, $\gamma_{2}=\gamma_{8}/\rbr{\beta \gamma_{6}}$, $\gamma_3 = 2 \beta \gamma_{5}$ and $\gamma_4 = 2 \beta \gamma_{6}$
we get the assertion. 
\end{proof}
\begin{corollary}Let $\alpha$ and $\gamma_{1},\ldots,\gamma_{4}$ be as
in the formulation of Theorem \ref{Tails}. If for $t>0$ we define
\begin{align*}
n_{1}(t) &:=\min\left\{ n:\frac{(\nu+n)^{2}}{2c^{2}}t\ge1,n=2,3,\ldots\right\} \ge 2, \\
n_{2}(t) &:=\min\left\{ n:\frac{(\nu+n)^{2}}{2c^{2}}t-\frac{(\nu+n_{1}(t))^{2}}{2c^{2}}t\ge1,n=2,3\ldots\right\} \ge n_{1}(t) +1 \ge 3
\end{align*}
and 
\[
F(t;\nu,c,z_0):=(1-\alpha)\left(\exp\rbr{-\frac{j_{\nu,1}^{2}}{2c^{2}}t}+\alpha\exp\rbr{-\frac{(\nu+2)^{2}}{2c^{2}}t}\min\rbr{\left\lceil \frac{1}{\ln(1/\alpha)}\right\rceil ,n_{2}(t)}\right)
\]
then  
\[
\P\rbr{\tau^{c}\ge\gamma_{1}t}\apprle F(t;\nu,c,z_0) \text { if } t\ge\gamma_{3}\E\tau^{c} \text{ and } 
\P\rbr{\tau^{c}\ge\gamma_{2}t}\gtrsim F(t;\nu,c,z_0) \text { if } t\ge\gamma_{4}\E\tau^{c}.
\]
\end{corollary}
\begin{proof}Let us notice that by the definition of $n_{2}(t)$ we have
\begin{equation}
\exp\rbr{-\frac{(\nu+n_{2}(t)-1)^2}{2c^{2}}t}\sim\exp\rbr{-\frac{(\nu+n_{1}(t))^2}{2c^{2}}t}\sim\exp\rbr{-\frac{(\nu+2)^2}{2c^{2}}t}.\label{eq:n_equiv}
\end{equation}
The second relation in \eqref{eq:n_equiv} follows from consideration of two cases: $n_1(t) = 2$ and $n_1(t) \ge 3$. In the latter case we have 
\[
\exp\rbr{-\frac{(\nu+2)^2}{2c^{2}}t} \sim \exp\rbr{-\frac{(\nu+n_1(t)-1)^2}{2c^{2}}t} \sim e^{-1} \sim 1, 
\]
but since $(\nu+n_1(t))^2 < 4 (\nu+n_1(t)-1)^2$ this also yields 
\[
 \exp\rbr{-\frac{(\nu+n_1(t))^2}{2c^{2}}t} \sim \exp\rbr{-\frac{(\nu+n_1(t)-1)^2}{2c^{2}}t} \sim 1.
\]
The inequalities $n_2(t) \ge 3$ and (\ref{eq:n_equiv})  also imply 
\begin{align*}
&\sum_{l=1}^{n_{2}(t)}\exp\rbr{-\frac{(\nu+l+1)^{2}}{2c^{2}}t} \le 3\sum_{l=1}^{n_{2}(t)-2}\exp\rbr{-\frac{(\nu+l+1)^{2}}{2c^{2}}t} \nonumber \\
& \sim \exp\rbr{-\frac{(\nu+2)^{2}}{2c^{2}}t} \rbr{n_2(t) -2} \sim \exp\rbr{-\frac{(\nu+2)^{2}}{2c^{2}}t} n_2(t)
\end{align*}
and on the other hand,
\begin{align*}
&\sum_{l=1}^{n_{2}(t)}\exp\rbr{-\frac{(\nu+l+1)^{2}}{2c^{2}}t} \ge \sum_{l=1}^{n_{2}(t)-1}\exp\rbr{-\frac{(\nu+l+1)^{2}}{2c^{2}}t} \nonumber \\
& \sim \exp\rbr{-\frac{(\nu+2)^{2}}{2c^{2}}t} \rbr{n_2(t) -1} \sim \exp\rbr{-\frac{(\nu+2)^{2}}{2c^{2}}t} n_2(t),
\end{align*}
so we have
\begin{equation} \label{lemacik}
\sum_{l=1}^{n_{2}(t)}\exp\rbr{-\frac{(\nu+l+1)^{2}}{2c^{2}}t} \sim \sum_{l=1}^{n_{2}(t)-2}\exp\rbr{-\frac{(\nu+l+1)^{2}}{2c^{2}}t} \sim \exp\rbr{-\frac{(\nu+2)^{2}}{2c^{2}}t} n_2(t).
\end{equation}
Let us also notice that if $l=1,2,\ldots,\left\lceil 1/\ln(1/\alpha)\right\rceil $
then $\alpha^{l-1}\sim1$ and that $\alpha^{\left\lceil 1/\ln(1/\alpha)\right\rceil }\le e^{-1}$.

If $\left\lceil 1/\ln(1/\alpha)\right\rceil \le n_{2}(t)-1$
then, using Theorem \ref{Tails}, we get that for $t\ge \gamma_{3}\E\tau^{c}$
\begin{align*}
&  \P\rbr{\tau^{c}\ge\gamma_{1}t} \apprle (1-\alpha)\exp\rbr{-\frac{j_{\nu,1}^{2}}{2c^{2}}t}+(1-\alpha)\alpha\sum_{n \ge 2}\exp\rbr{-\frac{(\nu+n)^{2}}{2c^{2}}t}\alpha^{n-2}\\
 & =(1-\alpha)\exp\rbr{-\frac{j_{\nu,n}^{2}}{2c^{2}}t}\\
 & \quad+(1-\alpha)\alpha\sum_{k=0}^{+\ns}\sum_{l=1}^{\left\lceil 1/\ln(1/\alpha)\right\rceil }\exp\rbr{-\frac{(\nu+k\left\lceil 1/\ln(1/\alpha)\right\rceil +l+1)^{2}}{2c^{2}}t}\alpha^{k\left\lceil 1/\ln(1/\alpha)\right\rceil +l-1}\\
 & \apprle (1-\alpha)\exp\rbr{-\frac{j_{\nu,n}^{2}}{2c^{2}}t}\\& \quad+(1-\alpha)\alpha\sum_{k=0}^{+\ns}\sum_{l=1}^{\left\lceil 1/\ln(1/\alpha)\right\rceil }\exp\rbr{-\frac{(\nu+k\left\lceil 1/\ln(1/\alpha)\right\rceil +l+1)^{2}}{2c^{2}}t}e^{-k}\\
 & \sim(1-\alpha)\exp\rbr{-\frac{j_{\nu,n}^{2}}{2c^{2}}t}+(1-\alpha)\alpha\sum_{l=1}^{\left\lceil 1/\ln(1/\alpha)\right\rceil }\exp\rbr{-\frac{(\nu+l+1)^{2}}{2c^{2}}t}\\
 & \sim(1-\alpha)\exp\rbr{-\frac{j_{\nu,n}^{2}}{2c^{2}}t}+(1-\alpha)\alpha\exp\rbr{-\frac{(\nu+2)^{2}}{2c^{2}}t}\left\lceil \frac{1}{\ln(1/\alpha)}\right\rceil,
\end{align*}
where in the last line we used \eqref{eq:n_equiv} (in the case $\left\lceil 1/\ln(1/\alpha)\right\rceil \le n_2(t) -2$) and \eqref{lemacik} (if $\left\lceil 1/\ln(1/\alpha)\right\rceil = n_2(t) -1$). Similar calculations yield $\P\rbr{\tau^{c}\ge\gamma_{2}t}\gtrsim F(t;\nu,c,\alpha)$
for $t\ge\gamma_{4}\E\tau^{c}$. 

Now we consider the case $\left\lceil 1/\ln({1}/{\alpha}) \right\rceil \ge n_{2}(t)$.
Since the function $n\mapsto(\nu+n)^{2}$ is convex for $n\ge2$,
for $k=1,2,\ldots$ and $l=1,2,\ldots,n_{2}(t)$ we have 
\begin{align*}
 & \frac{(\nu+n_{2}(t)\cdot k+l+1)^{2}}{2c^{2}}t-\frac{(\nu+n_{2}(t)\cdot(k-1)+l+1)^{2}}{2c^{2}}t\\
 & \ge\frac{(\nu+n_{2}(t))^{2}}{2c^{2}}t \ge \frac{(\nu+n_{2}(t))^{2}-(\nu+n_{1}(t))^{2}}{2c^{2}}t\ge1
\end{align*}
which implies 
\begin{equation}
\frac{(\nu+n_{2}(t)\cdot k+l+1)^{2}}{2c^{2}}t\ge\frac{(\nu+l+1)^{2}}{2c^{2}}t+k\label{eq:convexity_j}
\end{equation}
Now, using Theorem \ref{Tails},  (\ref{eq:convexity_j}) 
we get that for $t\ge\gamma_{3}\E\tau^{c}$
\begin{align*}
 & \P\rbr{\tau^{c}\ge\gamma_{1}t} \apprle (1-\alpha)\exp\rbr{-\frac{j_{\nu,n}^{2}}{2c^{2}}t}+(1-\alpha)\alpha\sum_{n \ge 2}\exp\rbr{-\frac{(\nu+n)^{2}}{2c^{2}}t}\alpha^{n-2}\\
 & = (1-\alpha)\exp\rbr{-\frac{j_{\nu,n}^{2}}{2c^{2}}t}+(1-\alpha)\alpha\sum_{k=0}^{+\ns}\sum_{l=1}^{n_{2}(t)}\exp\rbr{-\frac{(\nu+k\cdot n_{2}(t)+l+1)^{2}}{2c^{2}}t}\\
 & \sim(1-\alpha)\exp\rbr{-\frac{j_{\nu,n}^{2}}{2c^{2}}t}+(1-\alpha)\alpha\sum_{k=0}^{+\ns}\sum_{l=1}^{n_{2}(t)}\exp\rbr{-\frac{(\nu+l+1)^{2}}{2c^{2}}t}e^{-k}\\
 & \apprle (1-\alpha)\exp\rbr{-\frac{j_{\nu,n}^{2}}{2c^{2}}t}+(1-\alpha)\alpha\exp\rbr{-\frac{(\nu+2)^{2}}{2c^{2}}t}n_{2}(t),
\end{align*}
On the other hand, since $\left\lceil 1/\ln({1}/{\alpha})\right\rceil \ge n_{2}(t)$
then for $l=1,2,\ldots,n_{2}(t)$ we have $\alpha^{l-1}\sim1$ and
using Theorem \ref{Tails} and (\ref{eq:n_equiv}) we get that for
$t\ge\gamma_{4}\E\tau^{c}$
\begin{align*}
 & \P\rbr{\tau^{c}\ge\gamma_{2}t} \apprge (1-\alpha)\exp\rbr{-\frac{j_{\nu,n}^{2}}{2c^{2}}t}+(1-\alpha)\alpha\sum_{n \ge 2}\exp\rbr{-\frac{(\nu+n)^{2}}{2c^{2}}t}\alpha^{n-2}\\
 & \gtrsim(1-\alpha)\exp\rbr{-\frac{j_{\nu,n}^{2}}{2c^{2}}t}+(1-\alpha)\alpha\sum_{l=1}^{n_{2}(t)}\exp\rbr{-\frac{(\nu+l+1)^{2}}{2c^{2}}t}\alpha^{l-1}\\
 & \sim(1-\alpha)\exp\rbr{-\frac{j_{\nu,n}^{2}}{2c^{2}}t}+(1-\alpha)\alpha\exp\rbr{-\frac{(\nu+2)^{2}}{2c^{2}}t}n_{2}(t).
\end{align*}
\end{proof}
\begin{remark} If $n_{1}(t)\ge3$ and $\left\lceil 1/\ln(1/\alpha)\right\rceil <n_{2}(t)$
then 
\[
F(t;\nu,c,\alpha)\sim1.
\]

Indeed, if $n_{1}(t)\ge3$ then $(\nu+2)^{2}t/(2c^{2})<1$ and $\exp\rbr{-j_{\nu,1}^{2}t/(2c^{2})}\sim1$.
Thus, if $\left\lceil 1/\ln(1/\alpha)\right\rceil <n_{2}(t)$ then
\[
F(t;\nu,c,\alpha)\sim(1-\alpha)\left(1+\alpha\left\lceil \frac{1}{\ln({1}/{\alpha)}}\right\rceil \right)\sim(1-\alpha)\left(1+\frac{\alpha}{\ln({1}/{\alpha})}\right)\sim1.
\]
\end{remark}
\section{Estimtes of left tails of hitting times $\tau^{c}$ of Bessel processes}

In this section we present upper and lower bounds for left tails of
hitting times $\tau^{c}$. Similarly as in the previous section, they
will be obtained by application of different techniques. First, we
will apply exponential martingale technique to obtain upper bounds. Then, we will use Laplace transform estimate to obtain  lower bounds of the left tails for $\delta >1$. Finally we will utilize representation of $\tau^c$ as infinite convolution
of elementary mixtures of exponential distributions to obtain another upper and
lower bounds, valid for any $\delta >0$.

\subsection{Exponential martingale technique for upper bounds of the
left tails} \label{left_tail_up_b}

In this subsection we obtain upper bound for the left tail of hitting
times $\tau^{c}$. This way, for $d$-dimensional standard Brownian
motion we obtain estimates which have potential use in the Large Deviations
Theory. For example, \cite[Lemma 5.2.1]{Dembo:1998fu} provides some
estimate of $\P\rbr{\tau^{c/\sqrt{\varepsilon}}\le T}$ for fixed
$T$ and $c$ as $\varepsilon\ra0+$. In view of our results (see
Proposition \ref{Left_tails_rough}) and by scaling property of Brownian
motion we get optimal rate of this quantity for $\delta=2,3,\ldots$
as $\varepsilon\ra0+$. 

Let $X_{t}$, $Z_{t}$, $t\ge0$, be defined as in Section \ref{sec:Introduction}
and the process $Y_{t}$, $t\ge0$, be defined as in Subsection \ref{subsec:Exponential-martingale-technique_right_tails}.
Reasoning similarly as in Subsection \ref{subsec:Exponential-martingale-technique_right_tails} and using \eqref{eq:exp_estimate}
we get for any $u\in\R$ and $\lambda>0$ 
\begin{align*}
 & \P\rbr{c^{2}-\delta\cdot\tau^{c}-2\lambda c^{2}t-z_{0}\ge u\text{ and }\tau^{c}\le t}\\
 & =\mathbb{P}\left(\exp\left(\lambda\left(c^{2}-\delta\cdot\tau^{c}-2\lambda c^{2}t-z_{0}\right)\right)\ge e^{\lambda\cdot u}\text{ and }\tau^{c}\le t\right)\\
 & =\mathbb{P}\left(\exp\left(\lambda\left(Z_{t\wedge\tau^{c}}-\delta\cdot t\wedge\tau^{c}-2\lambda c^{2}t-z_{0}\right)\right)\ge e^{\lambda\cdot u}\text{ and }\tau^{c}\le t\right)\\
 & \le\E\exp\left(\lambda\left(Z_{t\wedge\tau^{c}}-\delta\cdot t\wedge\tau^{c}-2\lambda c^{2}t-z_{0}\right);\tau^{c}\le t\right)e^{-\lambda\cdot u}\\
 & \le\E\exp\left(\lambda\left(Z_{t\wedge\tau^{c}}-\delta\cdot t\wedge\tau^{c}-2\lambda c^{2}t-z_{0}\right)\right)e^{-\lambda\cdot u}
\le e^{-\lambda\cdot u}.
\end{align*}
Hence, for any $\eta\in(0,1)$, substituting $t=\left(1-\eta\right)\left(c^{2}-z_{0}\right)/\delta$,
$u=\eta\cdot\left(c^{2}-z_{0}\right)/2$, $\lambda=\delta\cdot\eta/(4c^{2}(1-\eta))$
we get
\[
\frac{c^{2}}{\delta}-2\lambda\frac{c^{2}}{\delta}t-\frac{u+z_{0}}{\delta}=\frac{c^{2}-z_{0}}{\delta}\left(1-\eta\right)
\]
and 
\begin{align*}
 & \P\rbr{c^{2}-\delta\cdot\tau^{c}-2\lambda c^{2}t-z_{0}\ge u\text{ and }\tau^{c}\le t}\\
 & =\P\rbr{\tau^{c}\le\frac{c^{2}}{\delta}-2\lambda\frac{c^{2}}{\delta}t-\frac{u+z_{0}}{\delta}\text{ and }\tau^{c}\le\frac{c^{2}-z_{0}}{\delta}\left(1-\eta\right)}\\
 & =\P\rbr{\tau^{c}\le\frac{c^{2}-z_{0}}{\delta}\left(1-\eta\right)}\\
 & \le\exp\rbr{-\frac{\delta\cdot\eta}{4c^{2}\left(1-\eta\right)}\frac{\eta\cdot\left(c^{2}-z_{0}\right)}{2}}=\exp\rbr{-\delta\frac{c^{2}-z_{0}}{c^{2}}\frac{\eta^{2}}{8(1-\eta)}}.
\end{align*}
Recalling that $\E\tau^{c}=\rbr{c^{2}-z_{0}}/\delta$ we can write
the just obtained estimate as 
\begin{align*}
\P\rbr{\tau^{c}\le\left(1-\eta\right)\E\tau^{c}} & \le\exp\rbr{-\delta\frac{c^{2}-z_{0}}{c^{2}}\frac{\eta^{2}}{8(1-\eta)}}.
\end{align*}

\subsection{Laplace transform technique for lower bounds of the
left tails}
\label{Laplace_left_t_l_b} 

To obtain estimates for lower bounds for the right tails we will use
another martingale. This technique will work for $\delta\ge2$ (equivalently
$\nu>0$) only. In the next subsection we will use similar technique
but we will apply directly the Laplace transform of $\tau^{c}$. For
$\delta>1$ the process $X$ satisfies the following SDE
\[
X_{t}=x_{0}+\beta_{t}+\frac{\delta-1}{2}\int_{0}^{t}X_{s}^{-1}\dd s,
\]
where $\beta_{t}$, $t\ge0$, is a standard Brownian motion, see \cite[Chap. XI, Exercise 1.26]{RevuzYor:2005}
(let us recall that $x_{0}=\sqrt{z_{0}}$). This implies that for
any $\lambda\in\R$, the process $\exp\rbr{\lambda\rbr{X_{t}-x_{0}-\frac{\delta-1}{2}\int_{0}^{t}X_{s}^{-1}\dd s}}=\exp\rbr{\lambda\beta_{t}}$
is a geometric Brownian motion and 
\[
U_{t}=\exp\rbr{\lambda\rbr{X_{t}-x_{0}-\frac{\delta-1}{2}\int_{0}^{t}X_{s}^{-1}\dd s}-\frac{1}{2}\lambda^{2}t}
\]
is a martingale. From this for any $\lambda>0$ and $t\ge0$ we get
\begin{align*}
1 & =U_{0}=\E U_{t\wedge\tau^{c}}=\E\exp\rbr{\lambda\rbr{X_{t\wedge\tau^{c}}-x_{0}-\frac{\delta-1}{2}\int_{0}^{t\wedge\tau^{c}}X_{s}^{-1}\dd s}-\frac{1}{2}\lambda^{2}t\wedge\tau^{c}}\\
 & \le\E\exp\rbr{\lambda\rbr{c-x_{0}-\frac{\delta-1}{2}c^{-1}t\wedge\tau^{c}}-\frac{1}{2}\lambda^{2}t\wedge\tau^{c}}
\end{align*}
which is equivalent to 
\[
\E\exp\rbr{-\rbr{\lambda\frac{\delta-1}{2c}+\frac{1}{2}\lambda^{2}}t\wedge\tau^{c}}\ge e^{-\lambda\rbr{c-x_{0}}}.
\]
Sending $t$ to $+\ns$ we get 
\begin{equation}
\E\exp\rbr{-\rbr{\lambda\frac{\delta-1}{2c}+\frac{1}{2}\lambda^{2}}\tau^{c}}\ge e^{-\lambda\rbr{c-x_{0}}}.\label{eq:estim_below_0}
\end{equation}
Denoting $u=\lambda\frac{\delta-1}{2c}+\frac{1}{2}\lambda^{2}>0$
we get 
\[
\lambda=\sqrt{\rbr{\frac{\delta-1}{2c}}^{2}+2u}-\frac{\delta-1}{2c}
\]
and 
\begin{equation}
\E\exp\rbr{-u\tau^{c}}\ge\exp\rbr{-\rbr{\sqrt{\rbr{\frac{\delta-1}{2c}}^{2}+2u}-\frac{\delta-1}{2c}}\rbr{c-x_{0}}}.\label{eq:estim_below}
\end{equation}
Notice that (\ref{eq:estim_below}) holds for any $u\ge0$. Next,
for any $t\ge0$ w estimate
\[
\E\exp\rbr{-u\tau^{c}}\le e^{-ut}\P\rbr{\tau^{c}>t}+\P\rbr{\tau^{c}\le t}=e^{-ut}+\rbr{1-e^{-ut}}\P\rbr{\tau^{c}\le t}
\]
which yields
\begin{equation}
\P\rbr{\tau^{c}\le t}\ge\frac{\E\exp\rbr{-u\tau^{c}}-e^{-ut}}{1-e^{-ut}}.\label{eq:estim_below1}
\end{equation}
If $A>0$ is some fixed number and $\lambda$ is the unique positive
solution of the equation
\begin{align}
u\cdot t=\rbr{\lambda\frac{\delta-1}{2c}+\frac{1}{2}\lambda^{2}}\cdot t & =A+\lambda\rbr{c-x_{0}}\label{eq:u}
\end{align}
then, by (\ref{eq:estim_below_0}) and (\ref{eq:estim_below1}) we
have
\begin{align}
\P\rbr{\tau^{c}\le t} & \ge\frac{e^{-\lambda\rbr{c-x_{0}}}\rbr{1-e^{-A}}}{1-e^{-A-\lambda\rbr{c-x_{0}}}} \ge e^{-\lambda\rbr{c-x_{0}}}\rbr{1-e^{-A}}.\label{eq:estim_below1-1}
\end{align}
Solving (\ref{eq:u}) for $\lambda$ we get 
\begin{align}
\lambda & =\frac{1}{t}\rbr{\sqrt{\rbr{\rbr{c-x_{0}}-\frac{\delta-1}{2c}t}^{2}+2At}+\rbr{c-x_{0}}-\frac{\delta-1}{2c}t}.\label{eq:lambda}
\end{align}
Substituting $t=2c\frac{c-x_{0}}{\delta-1}(1-\eta)\le 2 \frac{c^{2}-x_{0}^{2}}{\delta-1}(1-\eta)$
we get
\begin{align*}
\lambda & =\frac{\delta-1}{2c}\frac{\sqrt{\eta^{2}+4\frac{c}{c-x_{0}}\frac{A}{\delta-1}(1-\eta)}+\eta}{1-\eta}
\end{align*}
and this together with (\ref{eq:estim_below1-1}) yields
\begin{align}
 & \P\rbr{\tau^{c}\le2c\frac{c-x_{0}}{\delta-1}(1-\eta)} \nonumber \\
 & \ge\exp\rbr{-(\delta-1)\frac{c-x_{0}}{2c}\frac{\sqrt{\eta^{2}+4\frac{c}{c-x_{0}}\frac{A}{\delta-1}(1-\eta)}+\eta}{1-\eta}}\rbr{1-e^{-A}}. \label{CCC} 
\end{align}
By the inequality $\sqrt{a+b}\le\sqrt{a}+\sqrt{b}$ ($a,b\ge0$) we
have 
\begin{align}
 & (\delta-1)\frac{c-x_{0}}{2c}\frac{\sqrt{\eta^{2}+4\frac{c}{c-x_{0}}\frac{A}{\delta-1}(1-\eta)}+\eta}{1-\eta}  \le(\delta-1)\frac{c-x_{0}}{c}\frac{1}{1-\eta}+\sqrt{(\delta-1)\frac{c-x_{0}}{c}\frac{1}{1-\eta}}\sqrt{A}\label{eq:C}
\end{align}
and by the inequality $x+2\sqrt{x}\le2x+1$ ($x \ge 0$) we further have
\begin{equation}
(\delta-1)\frac{c-x_{0}}{c}\frac{1}{1-\eta}+\sqrt{(\delta-1)\frac{c-x_{0}}{c}\frac{1}{1-\eta}}\sqrt{4}\le2(\delta-1)\frac{c-x_{0}}{c}\frac{1}{1-\eta}+1.\label{eq:CC}
\end{equation}
By \eqref{eq:C} and \eqref{eq:CC}, substituting in \eqref{CCC} $A=4$ we have 
\begin{align*}
\P\rbr{\tau^{c}\le2c\frac{c-x_{0}}{\delta-1}(1-\eta)} \nonumber
 & \ge\exp\rbr{-2(\delta-1)\frac{c-x_{0}}{c}\frac{1}{1-\eta}-1}\rbr{1-e^{-4}} \nonumber \\
 & \ge \frac{1}{3} \exp\rbr{-2(\delta-1)\frac{c-x_{0}}{c} \frac{1}{1-\eta}}. 
 \end{align*}

To summarize the estimates obtained in this and the previous subsection
we state.

\begin{proposition} \label{Left_tails_rough} Let $c>0$ and $\tau^{c}$
be the hitting time of the $\delta$-dimensional Bessel process starting
from $x_{0}=\sqrt{z_{0}}\in[0,c)$, defined by (\ref{eq:tau_c_def}).
For any $\eta\in(0,1)$ and $\delta>0$ the following upper bound
holds: 
\[
\P\rbr{\tau^{c}\le\left(1-\eta\right)\E\tau^{c}} = \P\rbr{\tau^{c}\le\frac{c^{2}-x_{0}^{2}}{\delta}\left(1-\eta\right)}  \le\exp\rbr{-\frac{\delta}{8}\frac{c^{2}-x_{0}^{2}}{c^{2}}\frac{\eta^{2}}{1-\eta}}
\]
while for any $\eta\in(0,1)$ and $\delta>1$ the following lower
bounds holds:
\[
 \P\rbr{\tau^{c}\le2c\frac{c-x_{0}}{\delta-1}(1-\eta)} \ge \frac{1}{3} \exp\rbr{-2(\delta-1)\frac{c-x_{0}}{c} \frac{1}{1-\eta}}.
\]
\end{proposition} 

\subsection{Representation as infinite convolution
of elementary mixtures of exponential distributions and
bounds valid for any $\delta >0$.}
\label{Kent_left_t_l_b}  

In this section we apply representation of $\tau^c$ as infinite convolution
of elementary mixtures of exponential distributions to obtain
bounds valid for any $\delta >0$. To obtain these bounds we will also make use of the results of the previous two sections. 

In the case of small $\delta$, say $\delta\le2$ (equivalently $\nu=\delta/2-1\le0$)
we may reason in the following way. By \eqref{eq:law_tua_c} and \eqref{eq:law_tau_c1} we estimate 
\begin{equation}
\P\rbr{\tau^{c}\le s+u}\ge\P\rbr{\tau_{\nu,1}^{c}\le s}\P\rbr{\sum_{n \ge 2}\tau_{\nu,n}^{c}\le u}
\text{ and } \P\rbr{\tau_{\nu,1}^{c}\le s}=1-\rbr{1-\frac{z_{0}}{c^{2}}}\exp\rbr{-\frac{j_{\nu,1}^{2}}{2c^{2}}s}.\label{eq:B0}
\end{equation}

Let us notice that since the zeros of Bessel functions of the first
kind are interlacing, i.e. $j_{\nu,1}<j_{\nu+1,1}<j_{\nu,2}<j_{\nu+1,2}<\ldots$,
for any $u\ge0$ we naturally have
\[
\P\rbr{\tau_{\nu,n}^{c}\le u}\ge\P\rbr{\tau_{\nu+1,n-1}^{c}\le u}\text{ for }n=2,3,\ldots
\]
and hence 
\[
\P\rbr{\sum_{n \ge 2}\tau_{\nu,n}^{c}\le u}\ge\P\rbr{\sum_{n\ge 1}\tau_{\nu+1,n}^{c}\le u}.
\]
Let us denote $\tilde{\tau}^{c}=\sum_{n\ge 1}\tau_{\nu+1,n}^{c}$.
$\tilde{\tau}^{c}$ corresponds to the hitting time of a $2\rbr{\nu+2}=\delta+2$-dimensional
Bessel process and we have
\[
\E\tilde{\tau}^{c}=\frac{c^{2}-x_{0}^{2}}{\delta+2}.
\]
Thus, for $\eta\in(0,1)$ such that $1-\eta\ge\delta/(\delta+1)>\delta/(\delta+2)$
(from which follows $\frac{\delta+2}{\delta}(1-\eta)>1$) 
\begin{align}
 & \P\rbr{\sum_{n \ge 2}\tau_{\nu,n}^{c}\le4c\frac{c-x_{0}}{\delta}(1-\eta)}\ge\P\rbr{\tilde{\tau}^{c}\le4c\frac{c-x_{0}}{\delta}(1-\eta)}\nonumber \\
 & \ge\P\rbr{\tilde{\tau}^{c}\le2\rbr{c+x_{0}}\frac{c-x_{0}}{\delta}(1-\eta)}\nonumber \\
 & =\P\rbr{\tilde{\tau}^{c}\le2\frac{c^{2}-x_{0}^{2}}{\delta+2}\frac{\delta+2}{\delta}(1-\eta)}\ge\P\rbr{\tilde{\tau}^{c}<2\E\tilde{\tau}^{c}}\ge\frac{1}{2}.\label{eq:B}
\end{align}
Since $\delta+2>2$ we may apply Proposition \ref{Left_tails_rough} and for $1-\eta<\delta/(\delta+1)$
(which is equivalent to $\frac{\delta+1}{\delta}(1-\eta)<1$) obtain
\begin{align}
 & \P\rbr{\sum_{n \ge 2}\tau_{\nu,n}^{c}\le2c\frac{c-x_{0}}{\delta}(1-\eta)}\ge\P\rbr{\tilde{\tau}^{c}\le2c\frac{c-x_{0}}{\delta}(1-\eta)}\nonumber \\
 & =\P\rbr{\tilde{\tau}^{c}\le2c\frac{c-x_{0}}{\delta+1}\frac{\delta+1}{\delta}(1-\eta)}\nonumber \\
 & \ge\frac{1}{3}\exp\rbr{-2(\delta+1)\frac{c-x_{0}}{c}\frac{1}{\frac{\delta+1}{\delta}(1-\eta)}}=\frac{1}{3}\exp\rbr{-2\delta\frac{c-x_{0}}{c}\frac{1}{1-\eta}}.\label{eq:A}
\end{align}
By (\ref{eq:B}) and (\ref{eq:A}), for any $\eta\in\rbr{0,1}$ have
\begin{equation}
\P\rbr{\sum_{n \ge 2}\tau_{\nu,n}^{c}\le4c\frac{c-x_{0}}{\delta}(1-\eta)}\ge\frac{1}{3}\exp\rbr{-2\delta\frac{c-x_{0}}{c}\frac{1}{1-\eta}}.\label{eq:B1}
\end{equation}
From Proposition 3.1, (\ref{eq:B0}) and (\ref{eq:B1}), substituting in (\ref{eq:B0}) $s = \rbr{c^2-x_{0}^2}(1-\eta)/{\delta}$ (for upper bound) and $s = c\rbr{c-x_{0}}(1-\eta)/{\delta}$ (for lower bound), we obtain upper and lower bounds of left tails of $\tau^{c}$ such that these bounds as functions of $1-\eta$ are comparable when $1-\eta$ is rescaled by universal multiplicative constants. We have 

\begin{proposition} \label{prop:last} Let $c>0$ and $\tau^{c}$
be the hitting time of the $\delta$-dimensional Bessel process starting
from $x_{0}=\sqrt{z_{0}}\in[0,c)$, defined by (\ref{eq:tau_c_def}).
For any $\eta\in(0,1)$ and $\delta>0$ the following bounds hold
\begin{align*}
 & \P\rbr{\tau^{c}\le\frac{c^{2}-x_{0}^{2}}{\delta}\rbr{1-\eta}}\\
 & \le\min\cbr{\exp\rbr{-\delta\frac{c^{2}-x_{0}^{2}}{8c^{2}}\frac{\eta^{2}}{1-\eta}},1-\rbr{1-\frac{x_{0}^{2}}{c^{2}}}\exp\rbr{-\frac{j_{\nu,1}^{2}}{\delta}\frac{c^{2}-x_{0}^{2}}{2c^{2}}\rbr{1-\eta}}}
\end{align*}
and 
\begin{align*}
 & \P\rbr{\tau^{c}\le5c\frac{c-x_{0}}{\delta}(1-\eta)}\\
 & \ge\frac{1}{3}\exp\rbr{-2\delta\frac{c-x_{0}}{c}\frac{1}{1-\eta}}\cbr{1-\rbr{1-\frac{x_{0}^{2}}{c^{2}}}\exp\rbr{-\frac{j_{\nu,1}^{2}}{2\delta}\frac{c-x_{0}}{c}(1-\eta)}},
\end{align*}
where $\nu=\rbr{\delta/2}-1$ and $j_{\nu,1}$ is the first zero of
the Bessel function $J_{\nu}$ of the first kind.
\end{proposition}
\begin{remark}
For small $\delta$ and $x_{0}$ very close to $0$ the second
term on the right side of both estimates in Proposition \ref{prop:last} may be much smaller than
the first term. For example for $x_{0}=0$, recalling that (see \cite{Piessens:1984}) 
\[
\lim_{\delta\rightarrow0}\frac{j_{\nu,1}^{2}}{2\delta}=\lim_{\nu\rightarrow-1}\frac{j_{\nu,1}^{2}}{2\delta}=1
\]
for any fixed $\eta\in(0,1)$ we have 
\[
\lim_{\delta\rightarrow0}\frac{\exp\rbr{-\frac{\delta}{8}\frac{\eta^{2}}{1-\eta}}}{1-\exp\rbr{-\frac{j_{\nu,1}^{2}}{2\delta}\rbr{1-\eta}}}=\frac{1}{1-\exp\rbr{-\rbr{1-\eta}}}\ge\frac{1}{1-\eta}.
\]
\end{remark}
\subsection{Application to exit times of a $d$-dimensional standard Brownian motion}

By $B_{r}$ let us denote a closed ball in $\R^{d}$ centered at $0$:
$
B_{r}=\cbr{b\in\R^{d}:\left|b\right|\le r}.
$
Let $0<r<R$, $b_{0}\in B_{r}$ and $D$ be an open region in $\R^{d}$ such
that 
$
B_{r}\subseteq D\subseteq B_{R}.
$

Let $\tau^{D}$ denote exit time of $d$-dimensional standard Brownian
motion $B_t$, $t \ge 0$, such that $B_0 = b_0$, from the region $D$, that is
$
\tau^{D} = \inf\cbr{t>0: B_t \notin D}.
$ 
In this notation we have $\tau^{B_{r}}=\tau^{r}$
and $\tau^{B_{R}}=\tau^{R}$. 
 As a direct application of obtained estimates (Proposition \ref{Right_tails_rough} and Proposition \ref{Left_tails_rough}) we have the following corollary.
\begin{corollary} For any $\eta\in(0,1)$
\begin{align*}
\P\rbr{\tau^{D}\le\frac{r^{2}-\left|b_{0}\right|^{2}}{d}\left(1-\eta\right)}  \le\P\rbr{\tau^{B_{r}}\le\frac{r^{2}-\left|b_{0}\right|^{2}}{d}\left(1-\eta\right)}
 \le\exp\rbr{-d\frac{r^{2}-\left|b_{0}\right|^{2}}{r^{2}}\frac{\eta^{2}}{8(1-\eta)}}
\end{align*}
and for any $\eta>0$
\begin{align*}
\P\rbr{\tau^{D}\ge\frac{R^{2}-\left|b_{0}\right|^{2}}{d}\left(1+\eta\right)}  \le\P\rbr{\tau^{B_{R}}\ge\frac{R^{2}-\left|b_{0}\right|^{2}}{d}\left(1+\eta\right)}
  \le\exp\rbr{-d\frac{R^{2}-\left|b_{0}\right|^{2}}{R^{2}}\frac{\eta^{2}}{8(\eta+1)}}.
\end{align*}
\end{corollary}
\begin{proof} The estimates follow directly from Proposition \ref{Right_tails_rough} and Proposition \ref{Left_tails_rough} since (always)
$\tau^{D}\ge\tau^{B_{r}}$ and $\tau^{D}\le\tau^{B_{R}}$.
\end{proof}

{\bf Acknowledgments.} The research of both authors was funded by the National
Science Centre, Poland, under Grant No. 2019/35/B/ST1/042.


\begin{thebibliography}{KLRS07}

\bibitem[BR06]{Byczkowskietal:2006}
T.~Byczkowski and M.~Ryznar, \emph{Hitting distibution of geometric {B}rownian
  motion}, Stud. Math. \textbf{173} (2006), no.~1, 19--38.

\bibitem[Bre95]{Breen:1995}
Stephen Breen, \emph{Uniform upper and lower bounds on the zeros of {B}essel
  functions of the first kind}, J. Math. Anal. Appl. \textbf{196} (1995),
  no.~1, 1--17.

\bibitem[CT62]{CiesTaylor:1962}
Z.~Ciesielski and S.~J. Taylor, \emph{First passage times and sojourn times for
  brownian motion in space and the exact {H}ausdorff measure of the sample
  path}, Transactions of the American Mathematical Society \textbf{103} (1962),
  no.~3, 434--450.

\bibitem[dB87]{deBlassie1987}
Dante de~{Blassie}, \emph{{Stopping times of Bessel processes}}, {Ann. Probab.}
  \textbf{15} (1987), no.~3, 1044--1051 (English).

\bibitem[DZ98]{Dembo:1998fu}
A.~Dembo and O.~Zeitouni, \emph{{Large deviations techniques and
  applications}}, Springer, 1998.

\bibitem[GK95]{GluskinKwapien1995}
E.~D. {Gluskin} and S.~{Kwapie\'n}, \emph{{Tail and moment estimates for sums
  of independent random variables with logarithmically concave tails}}, {Stud.
  Math.} \textbf{114} (1995), no.~3, 303--309 (English).

\bibitem[HM12]{YamanaMatsumoto:2012}
Y.~Hamana and H.~Matsumoto, \emph{The probability densities of the first
  hitting times of bessel process}, J. of Math-for-Industry \textbf{4} (2012),
  91--95.

\bibitem[HM13]{YamanaMatsumoto:2013}
\bysame, \emph{The probability distributions of the first hitting times of
  {B}essel processes}, Trans. Amer. Math. Soc. \textbf{365} (2013), 5237--5257.

\bibitem[JS15]{jedidi:2015}
Wissem Jedidi and Thomas Simon, \emph{Diffusion hitting times and the
  bell-shape}, Statistics \& Probability Letters \textbf{102} (2015), 38 -- 41.

\bibitem[JW18]{Wisnie2018}
Jacek Jakubowski and Maciej Wi\'sniewolski, \emph{Invariance formulas for
  stopping times of squared bessel process}, Stochastic Analysis and
  Applications \textbf{36} (2018), no.~4, 671--699.

\bibitem[{Ken}80]{Kent:1980}
John~T. {Kent}, \emph{{Eigenvalue expansions for diffusion hitting times}}, {Z.
  Wahrscheinlichkeitstheor. Verw. Geb.} \textbf{52} (1980), 309--319 (English).

\bibitem[{Lat}97]{Latala1997}
Rafa{\l} {Lata{\l}a}, \emph{{Estimation of moments of sums of independent real
  random variables}}, {Ann. Probab.} \textbf{25} (1997), no.~3, 1502--1513
  (English).

\bibitem[L{\'e}v53]{Levy1953}
P.~L{\'e}vy, \emph{La mesure de {H}ausdorff de la courbe du mouvement
  brownien.}, Giorn. Ist. Ital. Attuari \textbf{16} (1953), 1--37.

\bibitem[MP10]{Peres:2010}
Peter M\"{o}rters and Yuval Peres, \emph{Brownian motion}, Cambridge University
  Press, Cambridge, 2010.

\bibitem[Pie84]{Piessens:1984}
R.~Piessens, \emph{A series expansion for the first positive zero of the
  {B}essel functions}, Mathematics of Computation \textbf{42} (1984), 195--197.

\bibitem[R{\"{o}}s90]{Roesler:1980}
U.~R{\"{o}}sler, \emph{Unimodality of passage times for one-dimensional strong
  {M}arkov processes}, Ann. Probab. \textbf{8} (1990), 143--172.

\bibitem[RY05]{RevuzYor:2005}
Daniel Revuz and Marc Yor, \emph{Continuous martingales and {B}rownian motion,
  3rd ed.}, Grundlehren der Mathematischen Wissenschaften, vol. 293,
  Springer-Verlag, Berlin, 2005. \MR{MR1083357 (92d:60053)}

\bibitem[Ser17]{Serafin:2017}
G.~Serafin, \emph{Exit times densities of the {B}essel process}, Proc. Amer.
  Math. Soc. \textbf{145(7)} (2017), 3165--3178.

\bibitem[TBR13]{Byczkowskietal:2012}
J.~Ma{\l}ecki T.~Byczkowski and M.~Ryznar, \emph{Hitting times of {B}essel
  processes}, Potential Analysis \textbf{38} (2013), 753?786.

\bibitem[TBS07]{Byczkowskietal:2007}
P.~Graczyk T.~Byczkowski and A.~St\'{o}s, \emph{Poisson kernels of half-spaces
  in real hyperbolic space}, Rev. Mat. Iberoam. \textbf{23} (2007), 85--126.

\bibitem[TGT19]{Grzywnyetal:2019}
M.~Ryznar T.~Grzywny and B.~Trojan, \emph{Asymptotic behaviour and estimates of
  slowly varying convolution semigroups}, Int. Math. Res. Not.
  \textbf{2019(23)} (2019), 7193--7258.

\bibitem[Wat44]{Watson:1944}
G.~N. Watson, \emph{{A Treatise on the Theory of {B}essel Functions}},
  Cambridge University Press, Cambridge, 1944.

\bibitem[Yam80]{Yamazato:1990}
M.~Yamazato, \emph{Hitting time distributions of single points for
  1-dimensional generalized diffusion processes}, Nagoya Math. J. \textbf{119}
  (1980), no.~4, 853--859.

\end{thebibliography}
\end{document}